\documentclass[11pt,reqno,twoside]{amsart}
\usepackage{amscd}
\usepackage{amsfonts}
\usepackage{amsmath}
\usepackage{amssymb}
\usepackage{amsthm}
\usepackage{fancyhdr}
\usepackage{latexsym}
\usepackage[colorlinks=true, pdfstartview=FitV, linkcolor=blue, citecolor=blue, urlcolor=blue]{hyperref}
\usepackage{enumitem}      
\usepackage{mathtools}            
\usepackage{upgreek}
\usepackage{caption} 
\usepackage{indentfirst} 
\usepackage{schemata} 
\usepackage{blindtext, rotating}   
\usepackage{soul} 
\usepackage{graphicx} 
\usepackage{float}
\setstcolor{red}
\usepackage{color}
\usepackage{tikz}
\usetikzlibrary{arrows.meta}
\usetikzlibrary{decorations.markings}
\tikzset{->-/.style={decoration={
  markings,
  mark=at position #1 with {\arrow{>}}},postaction={decorate}}}
  \tikzset{middlearrow/.style={
        decoration={markings,
            mark= at position 0.55 with {\arrow{#1}} ,
        },
        postaction={decorate}
    }
}          
\newcommand{\eee}[1]{\begin{equation}#1\end{equation}}
\newcommand{\sss}[1]{\begin{subequations}#1\end{subequations}}
\newcommand{\ddd}[1]{\begin{alignat}{2}#1\end{alignat}}

\newcommand{\nn}{\nonumber}
\newcommand{\p}{\partial}
\newcommand{\ve}{\varepsilon}

\definecolor{ddgreen}{RGB}{0,170,0}

\newcommand{\what}{\widehat}
\usepackage{scalerel,stackengine}
\stackMath
\newcommand\wwhat[1]{
\savestack{\tmpbox}{\stretchto{
  \scaleto{
    \scalerel*[\widthof{\ensuremath{#1}}]{\kern-.6pt\bigwedge\kern-.6pt}
    {\rule[-\textheight/2]{1ex}{\textheight}}
  }{\textheight} 
}{0.5ex}}
\stackon[1pt]{#1}{\tmpbox}
}
\makeatletter
\renewcommand\subsection{\@startsection{subsection}{2}%
  \z@{-0.8\linespacing\@plus-0.7\linespacing}{0.7\linespacing}%
  {\normalfont\bfseries}}
\makeatletter
\def\mathcolor#1#{\@mathcolor{#1}}
\def\@mathcolor#1#2#3{%
\protect\leavevmode
\begingroup
\color#1{#2}#3%
\endgroup
}
\makeatother
\theoremstyle{plain}

\theoremstyle{definition}

\newtheorem{remark}{Remark}[section]

\def\Xint#1{\mathchoice
   {\XXint\displaystyle\textstyle{#1}}%
   {\XXint\textstyle\scriptstyle{#1}}%
   {\XXint\scriptstyle\scriptscriptstyle{#1}}%
   {\XXint\scriptscriptstyle\scriptscriptstyle{#1}}%
   \!\int}
\def\XXint#1#2#3{{\setbox0=\hbox{$#1{#2#3}{\int}$}
     \vcenter{\hbox{$#2#3$}}\kern-.5\wd0}}

\def\dashint{\Xint-}

\DeclareMathSizes{12}{12}{7}{5}
\usepackage[framemethod=TikZ]{mdframed}
\usepackage[breakable, theorems, skins]{tcolorbox}
\tcbset{enhanced}

\numberwithin{figure}{section}
\numberwithin{equation}{section}
\usepackage{geometry}
\makeatletter
\def\l@section{\@tocline{1}{0pt}{1pc}{}{}}
\def\l@subsection{\@tocline{2}{0pt}{1pc}{4.6em}{}}
\def\l@subsubsection{\@tocline{3}{0pt}{1pc}{7.6em}{}}
\renewcommand{\tocsubsection}[3]{%
  \indentlabel{\@ifnotempty{#2}{\hspace*{2.3em}\makebox[2.3em][l]{%
    \ignorespaces#1 #2.\hfill}}}#3}
\renewcommand{\tocsubsubsection}[3]{%
  \indentlabel{\@ifnotempty{#2}{\hspace*{4.6em}\makebox[3em][l]{%
    \ignorespaces#1 #2.\hfill}}}#3}
\makeatother 
\usepackage{array,tabularx}
\newenvironment{conditions*}
  {\par\vspace{\abovedisplayskip}\noindent
   \tabularx{\columnwidth}{>{$}l<{$} @{\ : } >{\raggedright\arraybackslash}X}}
  {\endtabularx\par\vspace{\belowdisplayskip}}
\usepackage{mathrsfs}
\makeatletter
\@namedef{subjclassname@2020}{%
  \textup{2020} Mathematics Subject Classification}
\makeatother

\begin{document}

\title{The linearized classical Boussinesq system on the half-line}
\author{C. M. Johnston, Clarence T. Gartman \& Dionyssios Mantzavinos$^*$
\vskip 3mm
{\tiny Department of Mathematics, University of Kansas, Lawrence, KS 66045}}
\begin{abstract}
The linearization of the classical Boussinesq system is solved explicitly in the case of nonzero boundary conditions on the half-line. The analysis relies on the unified transform method of Fokas and is performed in two different frameworks: (i)  by exploiting the recently introduced extension of Fokas's method to systems of equations;  (ii) by expressing the linearized classical Boussinesq system as a single, higher-order equation which is then solved via the usual version of the unified transform. The resulting formula provides a novel representation for the solution of the linearized classical Boussinesq system on the half-line. Moreover, thanks to the uniform convergence at the boundary, the novel formula is shown to satisfy the linearized classical Boussinesq system as well as the prescribed initial and boundary data via a direct calculation.
\end{abstract}
\date{September 20, 2020. \textit{Revised}: November 8, 2020. $^*$\!\textit{Corresponding author}: mantzavinos@ku.edu}

\maketitle

\markboth{The linearized classical Boussinesq system on the half-line}{C. M. Johnston, Clarence T. Gartman \& Dionyssios Mantzavinos}


\section{Introduction}

The  classical  Boussinesq system
\begin{equation}\label{cbous}
\begin{array}{l}
v_t + u_x + \left(vu\right)_x = 0,
 \\
u_t + v_x + uu_x - u_{xxt} = 0,
\end{array}
\quad
\end{equation}
with $v=v(x, t)$,  $u=u(x,t)$, is a central model  in fluid dynamics that  captures the propagation of small-amplitude, weakly nonlinear shallow  waves on the free surface of an ideal irrotational fluid under the effect of gravity. The system \eqref{cbous} was first derived by Boussinesq \cite{b1872,b1877} as an asymptotic approximation (in the regime specified above) of the celebrated Euler equations of hydrodynamics.\footnote{More precisely, the system \eqref{cbous} was actually derived by Peregrine \cite{p1967}; Boussinesq had derived  a slightly modified version of that system which is not well-posed.} 
As such, it has been the subject of several works in the literature \cite{s1981,a1984,bcs2002,bcs2004,al2008,a2011,ad2012,l2013,mid2014,mtz2020,lw2020}  through various analytical as well as numerical techniques.  

The classical Boussinesq system \eqref{cbous} is nonlinear and hence challenging to study analytically. At the same time, the system is also dispersive. Hence, the existence and uniqueness of its solution (upon the prescription of suitable data)  can be established via the  powerful contraction mapping technique. A central role in the implementation of that technique is played by the solution map of the linear counterpart of the problem under consideration. More precisely, the (explicit) solution formula of the forced linear problem inspires an implicit mapping for the solution of the nonlinear problem; this mapping is then shown to be a contraction in an appropriate function space, thereby implying a unique solution (namely, the unique fixed point of the contraction) for the nonlinear problem. 
Therefore, the derivation of a linear solution formula which is effective for the purpose of function estimates is crucial in the investigation of the solvability of the nonlinear system \eqref{cbous}. 

Moreover, there exists a particular aspect of the  classical Boussinesq system \eqref{cbous} --- and of nonlinear dispersive systems in general --- which has not been explored much in the literature, namely, their formulation as   initial-boundary value problems (IBVPs) on domains that involve  a boundary, as opposed to the fully unbounded domain associated with the initial value problem. One of the main reasons behind the slow progress in the rigorous analysis of nonlinear dispersive IBVPs when compared to their associated initial value problems has been the absence of the Fourier transform in the IBVP setting. Indeed, we recall that in the case of the initial value problem  the linear solution formulae, which are essential for obtaining the basic estimates used in the contraction mapping technique, are easily derived by simply applying a Fourier transform in the spatial variable. Nevertheless, once a boundary is introduced in the problem (e.g. in the case of the half-line $\{x>0\}$) the spatial Fourier transform is no longer available and, even more, no classical transform exists that can produce linear solution formulae which are effective for the purpose of estimates. 

Motivated by the above, in this work we consider the linearization about zero of the classical Boussinesq system \eqref{cbous}, i.e. we set  $v(x, t) = \ve r(x, t)$ and $u(x, t) = \ve q(x, t)$  with $0<\ve\ll 1$, and study the resulting linearized classical Boussinesq system as an IBVP on the half-line with a Dirichlet boundary condition:
\sss{\label{lbous-ibvp}
\ddd{
&r_t + q_x = 0, \ \ q_t + r_x - q_{xxt} = 0,  &&x>0, \ t>0,
\label{lbous-sys}
\\
&r(x, 0) = r_0(x), \ \ q(x, 0) = q_0(x), \quad &&x\geqslant 0,
\label{lbous-ic}
\\
&q(0, t) = g_0(t), \ \ t\geqslant 0,
\label{lbous-bc}
}
}
where, for the purpose of this work, we assume  initial and boundary data with  sufficient smoothness and decay at infinity (e.g. in the Schwartz class).\footnote{There exist works in the literature that study dispersive IBVPs in the case of rough data. This task, however, lies beyond the scope of the present article.} In particular, we assume  compatibility of the initial and boundary data at the origin, namely $r_0(0) = g_0(0)$. 
We remark that the above IBVP is supplemented with just one boundary condition for $q$ and no boundary condition for $r$. Although it is not a priori clear that this choice of data is admissible (and sufficient), our analysis will reveal that this is indeed the case.

As noted earlier, no classical spatial transform can produce a solution formula for IBVP \eqref{lbous-ibvp} which is appropriate for analyzing the corresponding nonlinear IBVP. We emphasize that this is not a pathogeny of the linearized classical Boussinesq system, but rather a challenge which is present across the whole spectrum of linear dispersive IBVPs with nonzero boundary conditions. The resolution to this important obstacle in the study of linear (and, eventually, nonlinear) dispersive IBVPs was provided in 1997 by Fokas  \cite{f1997}, who introduced the now well-established unified transform method (UTM),  also known in the literature as the Fokas method. This novel method essentially provides the direct analogue of the Fourier transform in the IBVP framework. It relies on exploiting certain symmetries of the dispersion relation of the problem together with the ability to deform certain paths of integration to appropriate contours in the complex spectral plane. Over the last twenty years, UTM has been widely used for a plethora of linear as well as nonlinear evolution and elliptic equations, formulated on various domains in one or higher dimensions, and with a broad range of (admissible) boundary conditions --- see, for example, the research articles \cite{fk2003,fi2004,af2005,ff2008,ffss2009,ssf2010,fl2012,hm2015a,dss2016,cff2018,ko2020}, the books \cite{fbook,fp2015} and the review articles \cite{fs2012,dtv2014}.  

Fairly recently, Deconinck, Guo, Shlizerman and Vasan extended the \textit{linear} component of UTM   to systems of equations \cite{dgsv2018}. In this work, we shall exploit that recent progress   in order to derive the novel, explicit representation \eqref{lcbous-sol-fin} for the solution of the linearized classical Boussinesq IBVP \eqref{lbous-ibvp}. Our derivation will be done under the assumption of existence of solution. Nevertheless, taking advantage of one of the key features of UTM, namely the uniform convergence of its solution formulae at the boundary, we shall explicitly demonstrate (via a direct calculation) that our novel formula does indeed satisfy IBVP \eqref{lbous-ibvp} and is, in fact, a classical solution. Furthermore, we shall also provide an alternative way of solving IBVP \eqref{lbous-ibvp} by converting the linearized classical Boussinesq system into a single equation. A slight downside of the latter approach is that the resulting single equation involves a second-order time derivative, which makes the analysis somewhat more tedious. At the same time, the latter approach illustrates that UTM as a method is equally effective in both the single equation and the system frameworks.

It should be noted that IBVP \eqref{lbous-ibvp} has previously been considered  by Fokas and Pelloni in \cite{fp2005}. However, the method employed in that paper was the \textit{nonlinear} component of UTM, which relies on expressing the linearized classical Boussinesq system as a Lax pair and then using ideas inspired by the  inverse scattering transform in order to associate the solution of IBVP \eqref{lbous-ibvp} to that of a (scalar) Riemann-Hilbert problem. This Riemann-Hilbert problem is then solved explicitly with the help of Plemelj's formulae to yield the solution to problem \eqref{lbous-ibvp}. Nevertheless, the highly technical aspects of the approach of \cite{fp2005} make it less accessible to the broader applied sciences community. Furthermore, the solution representation produced in \cite{fp2005} involves certain principal value integrals, which arise as byproducts of the Plemelj formulae. This feature does not seem  convenient regarding  (i) the derivation of linear estimates for the contraction mapping analysis of the nonlinear system \eqref{cbous}, and (ii) numerical considerations. On the contrary, the new solution representation derived in the present work does not involve principal value integrals and, more importantly, relies solely on the \textit{linear} component of UTM, which only requires knowledge of the Fourier transform and of Cauchy's theorem from complex analysis.

\vskip 3mm
\noindent
\textbf{Organization of the article.} In Section \ref{gr-s}, starting from the linearized classical Boussinesq  IBVP \eqref{lbous-ibvp} we derive an important spectral identity which plays a central role in UTM and is known as the global relation. In Section \ref{sol-s}, we combine the global relation with appropriate deformations in the complex spectral plane in order to obtain the explicit solution formula of problem \eqref{lbous-ibvp}. Then, in Section \ref{one-eq-s}, we revisit the problem by converting the linearized classical Boussinesq system into a single equation which we then solve via the standard version of UTM. This offers a first, indirect way of corroborating our novel solution formula. The direct, explicit verification of our  formula is then presented in detail in Section \ref{ver-s}. Finally, some concluding remarks are provided in Section \ref{conc-s}.

\section{Derivation of the global relation}
\label{gr-s}

We define the half-line Fourier transform  by
\sss{
\eee{\label{ft-def}
\hat f(k)  =\int_{x=0}^{\infty}e^{-ikx}  f(x)dx, \quad \text{Im}(k) \leqslant 0,
}
with inverse
\eee{
f(x)  = \frac{1}{2\pi} \int_{k\in \mathbb R} e^{ikx}  \hat f(k) dk, \quad  x\geqslant 0.
\label{ift-def}
} 
}
We note that, unlike the standard Fourier transform, which is only valid on the real line, the half-line Fourier transform \eqref{ft-def} is valid on the closure of the lower half of the complex $k$-plane due to the fact that $x\geqslant 0$.
Then, applying \eqref{ft-def} to the linearized classical Boussinesq system \eqref{lbous-sys} we obtain
\begin{equation}\label{ode-qr0}
\arraycolsep=1.6pt
\def\arraystretch{1.2}
\begin{array}{rcl}
\p_t \, \hat r(k, t) &=& -ik\, \hat q(k, t)+ g_0(t), 
\\
\left(1+k^2\right) \p_t \, \hat q(k, t) &=& -ik \, \hat r(k, t) + h_0(t) -  g_1'(t) - ik g_0'(t),
\end{array}
\end{equation}
where $'$ denotes differentiation with respect to $t$ and we have introduced the   notation
\begin{equation}
g_j(t) =\partial_x^{j} q(0,t), \quad h_j(t) =\partial_x^{j} r(0,t), \quad j \in \mathbb N\cup \{0\}.
\end{equation}
For $k\neq \pm i$, we may express the system of ordinary differential equations \eqref{ode-qr0} in matrix form as
\begin{equation}\label{ode-q-vec0}
\p_t \, \hat{\mathbf q}(k, t) = A(k) \, \hat{\mathbf q}(k, t)+\mathbf{g}(k, t),
\end{equation}
where
\ddd{
&\hat{\mathbf q}(k, t) = \arraycolsep=3pt \left(\begin{array}{c} \hat r(k, t) \\ \hat q(k, t) \end{array}\right),
\quad
A(k) = -ik 
\left(\begin{array}{cc}
0 & 1 \\
\frac{1}{1+k^2} & 0 
\end{array}\right),
\nn\\
&
\mathbf g(k, t) = 
\def\arraystretch{1.2}
\left(\begin{array}{c}
g_0(t) \\
\frac{1}{1+k^2} \left[h_0(t) - g_1'(t) - ik g_0'(t)\right]
\end{array}\right).
}
Next, recalling the definition of the matrix exponential  $e^A := \sum_{j = 0}^{\infty}\frac{A^j}{j!}$, we integrate  \eqref{ode-q-vec0} to obtain the following spectral identity which in the UTM terminology is known as the \textit{global relation} (since it only involves  integrals of the vector $\mathbf q = (r, q)^T$ and its initial and boundary values):
\begin{equation}\label{l-serre-gr-vec0}
\widehat{\mathbf q}(k, t)=e^{A(k) t}\, \widehat{\mathbf q}_0(k) + \int_{\tau=0}^{t}e^{A(k)(t-\tau)}\, \mathbf{g}(k, \tau)\, d\tau, 
\quad \text{Im}(k) \leqslant 0, \  k\neq - i.
\end{equation}
We emphasize that the global relation  is valid for all $\text{Im}(k) \leqslant 0$ with $k\neq -i$ because the half-line Fourier transform \eqref{ft-def} makes sense for all  $\text{Im}(k) \leqslant 0$.

\begin{figure}[ht]
\centering
\begin{tikzpicture}[scale=1.8, rotate=0]
\draw [->,>=Stealth] (-0.75, 0) -- (1.5, 0);
\draw [->,>=Stealth]  (0, -1.25) -- (0, 1.25);
\draw[line width=0.6mm] (0, -0.75) -- (0, 0.75);
\draw[dashed, rotate around={-19:(0,0.75)}] (0, 0.75) -- (1.05,0.75);
\draw[dashed, rotate around={-39:(0,-0.75)}] (0, -0.75) -- (0,0.75);
\draw [domain=-90:-19, ->, >=stealth] plot ({0.25*cos(\x)}, {0.75+0.25*sin(\x)});
\node[] at (0.2, 0.45) {\fontsize{9}{9} $\theta_2$};
\draw [domain=-90:49, ->, >=stealth] plot ({0.2*cos(\x)}, {-0.75+0.2*sin(\x)});
\node[] at (0.28, -0.87) {\fontsize{9}{9} $\theta_1$};
\node[] at (1.05, 0.45) {\fontsize{9}{9} $k$};
\node[] at (-0.17, 0.75) {\fontsize{11}{11} $i$};
\node[] at (-0.25, -0.73) {\fontsize{11}{11} $-i$};
\filldraw (0,0.75) circle (0.75pt);
\filldraw (0,-0.75) circle (0.75pt);
\filldraw (0.95, 0.42) circle (0.75pt);
\node[] at (-0.17, 0.15) {\fontsize{9}{9} $\mathcal B$};
\end{tikzpicture}
\caption{The definition \eqref{sqrt-def} of the complex square root   $\mu(k) = \left(1+k^2\right)^{\frac 12}$ as a single-valued function by taking a branch cut along the segment $\mathcal B = i \left[-1, 1\right]$.}
\label{branch-cut-f}
\end{figure}
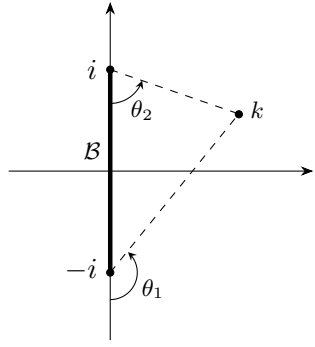

It turns out convenient to express the global relation \eqref{l-serre-gr-vec0}  in component form. For this purpose, we first diagonalize the matrix $A$ as
\eee{
A = P D P^{-1}
}
with
\eee{\label{om-def-eq}
P = \left(\begin{array}{cc}   1   & 1 \\ -\frac{1}{\mu} &  \frac{1}{\mu}\end{array}\right),
\quad
D = \left(\begin{array}{cc}   i\omega   & 0 \\ 0 &  -i\omega  \end{array}\right),
\quad
 \omega = \omega(k) := \frac{k}{\mu(k)},
}
where the complex square root
\eee{
\mu(k) := \left(1+k^2\right)^{\frac 12}
}
is made single-valued by taking a branch along the segment $\mathcal B:= i[-1, 1]$. In particular, we define
\eee{\label{sqrt-def}
\left(1+k^2\right)^{\frac 12}
= \sqrt{\left|1+k^2\right|} \, e^{i\left(\theta_1+\theta_2-\pi\right)/2}, \quad k\notin \mathcal B,
}
with the angles $\theta_1, \theta_2 \in [0, 2\pi)$ as shown in Figure \ref{branch-cut-f} and for $k\in \mathcal B$ we identify $\mu(k)$ by its limit  from the right. 

The above diagonalization allows us to write $e^{At}$ explicitly as a matrix:
\eee{\label{eAt0}
e^{At}
=
P e^{Dt} P^{-1}
=
\frac 12 
\def\arraystretch{1.2}
\left(\begin{array}{cc}  e^{i\omega t} + e^{-i\omega  t}  & -  \mu  \left(e^{i\omega  t}-e^{-i\omega t}\right) 
\\
-\frac{1}{\mu} \left(e^{i\omega t}-e^{-i\omega  t}\right) &  e^{i\omega  t}+e^{-i\omega t} 
\end{array}\right).
}
In turn, for $\text{Im}(k) \leqslant 0$ and $k\neq -i$ the global relation \eqref{l-serre-gr-vec0} can be expressed in component form as
\sss{\label{serre-gr-U=0-0}
\ddd{
\hat r(k, t) 
&= \frac 12 \left[ \left(e^{i\omega t}+e^{-i\omega t}\right) \hat r_0(k) -   \mu(k)  \left(e^{i\omega t}-e^{-i\omega t}\right) \hat q_0(k) \right]
\\
&\quad
+
\frac{1}{2\mu(k)} \, \bigg\{ -\left[ e^{i\omega t} \,  \widetilde h_0(\omega, t) - e^{-i\omega t} \, \widetilde h_0(-\omega, t)\right]
+
\mu(k) \left[ e^{i\omega t} \,  \widetilde g_0(\omega, t) + e^{-i\omega t} \,  \widetilde g_0(-\omega, t) \right]
\nn\\
&\hskip 2.25cm
+
ik \left[e^{i\omega t} \,  \widetilde{g_0'}(\omega, t) - e^{-i\omega t} \,  \widetilde{g_0'}(-\omega, t)\right]
+
 \left[e^{i\omega t} \,  \widetilde{g_1'}(\omega, t) - e^{-i\omega t} \,  \widetilde{g_1'}(-\omega, t)\right] \bigg\}
 \nn
}
and
\ddd{
\hat q(k, t) 
&= \frac 12 \left[ \left(e^{i\omega t}+e^{-i\omega t}\right) \hat q_0(k) - \frac{1}{\mu(k)}  \left(e^{i\omega t}-e^{-i\omega t}\right) \hat r_0(k) \right]
\\
&\quad
+ \frac{1}{2\left(1+k^2\right)} \, \bigg\{ \left[ e^{i\omega t} \, \widetilde h_0(\omega, t) + e^{-i\omega t} \, \widetilde h_0(-\omega, t) \right]
-  \mu(k) \left[ e^{i\omega t} \, \widetilde g_0(\omega, t) - e^{-i\omega t} \, \widetilde g_0(-\omega, t) \right]
\nn\\
&\hskip 2.85cm
- ik \left[e^{i\omega t} \,  \widetilde{g_0'}(\omega, t) + e^{-i\omega t} \,  \widetilde{g_0'}(-\omega, t)\right]
- \left[e^{i\omega t} \,  \widetilde{g_1'}(\omega, t) + e^{-i\omega t} \,  \widetilde{g_1'}(-\omega, t)\right] \bigg\},
\nn
}
}
where we have introduced the notation
\eee{
\widetilde f(\omega, t) = \int_{\tau=0}^t e^{-i\omega\tau} f(\tau) d\tau.
}

\section{Elimination of the unknown boundary values}
\label{sol-s}

The global relations \eqref{serre-gr-U=0-0} involve three boundary values, one for $r$ and two for $q$. As usual in the context of UTM,  escaping to the complex $k$-plane by means of Cauchy's theorem will allow us to eliminate two of those boundary values and thereby derive an effective solution representation, in the sense that it will only involve the boundary value $q(0, t) = g_0(t)$  prescribed as a datum in IBVP \eqref{lbous-ibvp}. This elimination procedure illustrates the ability of UTM to indicate which boundary values are admissible as data for a well-posed problem.

We begin by observing that the first of equations \eqref{lbous-ibvp} evaluated at $x=0$ implies
\eee{\label{g1-h0-rel}
g_1(t) = - h_0'(t).
} 
The above formal evaluation is done (as noted earlier) under the assumption of existence of a smooth solution; it can be verified a posteriori by direct evaluation of formula \eqref{lcbous-sol-fin} using the methods of Section \ref{ver-s}.
Hence, the boundary value $g_1$ can be eliminated from the global relations \eqref{serre-gr-U=0-0}, which now read
\sss{\label{serre-gr-U=0}
\ddd{
\hat r(k, t) 
&= \frac 12 \left[ \left(e^{i\omega t}+e^{-i\omega t}\right) \hat r_0(k) -    \mu(k)  \left(e^{i\omega t}-e^{-i\omega t}\right) \hat q_0(k) \right]
\\
&\quad
+
\frac{1}{2 \mu(k)} \, \bigg\{ -\left[ e^{i\omega t} \,  \widetilde h_0(\omega, t) - e^{-i\omega t} \, \widetilde h_0(-\omega, t)\right]
-
\left[e^{i\omega t} \,  \widetilde{h_0''}(\omega, t) - e^{-i\omega t} \,  \widetilde{h_0''}(-\omega, t)\right]
\nn\\
&\hskip 2.25cm
+
\mu(k)  \left[ e^{i\omega t} \,  \widetilde g_0(\omega, t) + e^{-i\omega t} \,  \widetilde g_0(-\omega, t) \right]
+ ik \left[e^{i\omega t} \,  \widetilde{g_0'}(\omega, t) - e^{-i\omega t} \,  \widetilde{g_0'}(-\omega, t)\right] \bigg\}
\nn
}
and
\ddd{\label{serre-gr-U=0-r}
\hat q(k, t) 
&= \frac 12 \left[ \left(e^{i\omega t}+e^{-i\omega t}\right) \hat q_0(k) - \frac{1}{\mu(k)}  \left(e^{i\omega t}-e^{-i\omega t}\right) \hat r_0(k) \right]
\\
&\quad
+ \frac{1}{2\left(1+k^2\right)} \, \bigg\{ \left[ e^{i\omega t} \, \widetilde h_0(\omega, t) + e^{-i\omega t} \, \widetilde h_0(-\omega, t) \right]
+ \left[e^{i\omega t} \,  \widetilde{h_0''}(\omega, t) + e^{-i\omega t} \,  \widetilde{h_0''}(-\omega, t)\right]
\nn\\
&\hskip 2.5cm
- \mu(k) \left[ e^{i\omega t} \, \widetilde g_0(\omega, t) - e^{-i\omega t} \, \widetilde g_0(-\omega, t) \right]
- ik  \left[e^{i\omega t} \,  \widetilde{g_0'}(\omega, t) + e^{-i\omega t} \,  \widetilde{g_0'}(-\omega, t)\right] \bigg\}.
\nn
}
}

Recall  that the global relations \eqref{serre-gr-U=0} are valid for  $\text{Im}(k)\leqslant 0$ with $k\neq -i$. Thus, using these expressions for $k\in \mathbb R$ in the  Fourier inversion \eqref{ift-def} yields the following integral representations  for the two components of system \eqref{lbous-ibvp}:
\sss{\label{serre-ir-U=0}
\ddd{
r(x, t) 
&= \frac{1}{2\pi}  \int_{k\in \mathbb R} e^{ikx} \, \frac 12 \left[ \left(e^{i\omega t}+e^{-i\omega t}\right) \hat r_0(k) - \mu(k)  \left(e^{i\omega t}-e^{-i\omega t}\right) \hat q_0(k) \right]  dk
\nn\\
&\quad
+ \frac{1}{2\pi} \int_{k\in \mathbb R} e^{ikx} \, \frac{1}{2\mu(k)} \, \bigg\{ - \left[e^{i\omega t} \,  \widetilde h_0(\omega, t) - e^{-i\omega t} \, \widetilde h_0(-\omega, t)\right] 
-  \left[e^{i\omega t} \,  \widetilde{h_0''}(\omega, t) - e^{-i\omega t} \,  \widetilde{h_0''}(-\omega, t)\right]  
\nn\\
&\quad
+  \mu(k) \left[ e^{i\omega t} \,  \widetilde g_0(\omega, t) + e^{-i\omega t} \,  \widetilde g_0(-\omega, t) \right] 
+  ik   \left[e^{i\omega t} \,  \widetilde{g_0'}(\omega, t) - e^{-i\omega t} \,  \widetilde{g_0'}(-\omega, t)\right] \bigg\} \, dk
}
and
\ddd{
&
q(x, t) 
= \frac{1}{2\pi} \int_{k\in \mathbb R} e^{ikx} \, \frac 12 \left[ \left(e^{i\omega t}+e^{-i\omega t}\right) \hat q_0(k) - \frac{1}{\mu(k)}  \left(e^{i\omega t}-e^{-i\omega t}\right) \hat r_0(k) \right] dk
\nn\\
&
+ \frac{1}{2\pi}   \int_{k\in \mathbb R} e^{ikx} \, \frac{1}{2\left(1+k^2\right)} \, \bigg\{  \left[ e^{i\omega t} \, \widetilde h_0(\omega, t) + e^{-i\omega t} \, \widetilde h_0(-\omega, t) \right] 
+ \left[e^{i\omega t} \,  \widetilde{h_0''}(\omega, t) + e^{-i\omega t} \,  \widetilde{h_0''}(-\omega, t)\right] \bigg\}\, dk
\nn\\
&
- \mu(k) \left[ e^{i\omega t} \, \widetilde g_0(\omega, t) - e^{-i\omega t} \, \widetilde g_0(-\omega, t) \right] 
- ik \left[e^{i\omega t} \,  \widetilde{g_0'}(\omega, t) + e^{-i\omega t} \,  \widetilde{g_0'}(-\omega, t)\right] \bigg\}\, dk.
}
}

\begin{remark}[Crossing the branch cut]\label{bcut-r}
Although  $\omega$ inherits the branch cut $\mathcal B$ from $\mu$, the integrands involved in  \eqref{serre-ir-U=0} are entire  in $k$. Indeed, denoting by $\mu^+$ and $\mu^-$ the limits of $\mu$ as $k$ approaches $\mathcal B$ from the left and from the right respectively,   according to \eqref{sqrt-def} we have $\mu^+  =  -\mu^-$, 
i.e. $\mu$ changes sign across $\mathcal B$. In turn, $\omega^+  = -\omega^-$ and hence the functions
$e^{i\omega t}  + e^{-i\omega t}$, $\frac{e^{i\omega t} - e^{-i\omega t}}{\mu}$, $\mu \left(e^{i\omega t} - e^{-i\omega t}\right)$ and $\mu^2$
are continuous across $\mathcal B$ and, therefore, the paths of integration in \eqref{serre-ir-U=0} are allowed to cross $\mathcal B$.
\end{remark}

The integral representations \eqref{serre-ir-U=0} involve two boundary values, $q(0, t) = g_0(t)$ and $r(0, t) = h_0(t)$. However, only the first one is prescribed as a boundary condition in problem \eqref{lbous-ibvp}, which is the reason why \eqref{serre-ir-U=0} is not an effective solution formula. Next, we shall eliminate from \eqref{serre-ir-U=0} the unknown boundary value $h_0(t)$ or, more precisely, the transforms $\widetilde h_0$ and $\widetilde{h_0''}$,   by using a symmetry of the global relations \eqref{serre-gr-U=0}. 
More specifically, we note that the transformation $k \mapsto -k$ leaves $\omega$ invariant.  This is because the definition \eqref{sqrt-def} implies\footnote{An easy way to see that $\mu(-k) = -\mu(k)$ is to observe that the branch cut for $\mu$ is such that $\mu(k) \simeq k$ as $|k|\to \infty$.} $\mu(-k) = -\mu(k)$ and hence $\omega(-k) = \omega(k)$.  Thus, under the transformation $k \mapsto -k$ the global relations \eqref{serre-gr-U=0} yield  the identities
\sss{\label{serre-gr-U=0-sym}
\ddd{\label{serre-gr-U=0-sym-a}
\hat r(-k, t) 
&= \frac 12 \left[ \left(e^{i\omega t}+e^{-i\omega t}\right) \hat r_0(-k) +  \mu(k)  \left(e^{i\omega t}-e^{-i\omega t}\right) \hat q_0(-k) \right]
\\
&\quad
+
\frac{1}{2\mu(k)} \, \bigg\{ \left[ e^{i\omega t} \,  \widetilde h_0(\omega, t) - e^{-i\omega t} \, \widetilde h_0(-\omega, t)\right]
+ \left[e^{i\omega t} \,  \widetilde{h_0''}(\omega, t) - e^{-i\omega t} \,  \widetilde{h_0''}(-\omega, t)\right]
\nn\\
&\hskip 2.2cm
+
\mu(k) \left[ e^{i\omega t} \,  \widetilde g_0(\omega, t) + e^{-i\omega t} \,  \widetilde g_0(-\omega, t) \right]
+
ik \left[e^{i\omega t} \,  \widetilde{g_0'}(\omega, t) - e^{-i\omega t} \,  \widetilde{g_0'}(-\omega, t)\right] \bigg\}
\nn
}
and
\ddd{\label{serre-gr-U=0-sym-b}
\hat q(-k, t) 
&= \frac 12 \left[ \left(e^{i\omega t}+e^{-i\omega t}\right) \hat q_0(-k) + \frac{1}{\mu(k)}  \left(e^{i\omega t}-e^{-i\omega t}\right) \hat r_0(-k) \right]
\\
&\quad
+ \frac{1}{2\left(1+k^2\right)} \, \bigg\{ \left[ e^{i\omega t} \, \widetilde h_0(\omega, t) + e^{-i\omega t} \, \widetilde h_0(-\omega, t)\right]
+  \left[e^{i\omega t} \,  \widetilde{h_0''}(\omega, t) + e^{-i\omega t} \,  \widetilde{h_0''}(-\omega, t)\right]
\nn\\
&\hskip 2.3cm
+ \mu(k) \left[ e^{i\omega t} \, \widetilde g_0(\omega, t) - e^{-i\omega t} \, \widetilde g_0(-\omega, t)\right]
+ ik \left[e^{i\omega t} \,  \widetilde{g_0'}(\omega, t) + e^{-i\omega t} \,  \widetilde{g_0'}(-\omega, t)\right] \bigg\}.
\nn
}
}
Since the global relations \eqref{serre-gr-U=0} are valid for $\text{Im}(k) \leqslant 0$ with $k\neq -i$, the identities \eqref{serre-gr-U=0-sym} hold for $\text{Im}(k) \geqslant 0$ with $k\neq i$. Thus, they can be readily employed for eliminating the unknown transforms $\widetilde h_0$ and $\widetilde{h_0''}$. 
However, before doing so, it turns out useful to first deform the contours of integration of the integrals  in \eqref{serre-ir-U=0} which involve the boundary values from the real axis to the complex $k$-plane.

\begin{figure}[ht]
\centering
\begin{tikzpicture}[scale=2, rotate=0]
\draw [->,>=Stealth] (-0.75, 0.205) -- (1.2, 0.205);
\draw [->,>=Stealth]  (0, -0.65) -- (0, 1.25);
\draw[line width=0.6mm] (0, -0.3) -- (0, 0.75);
 \draw  plot[smooth, tension=.7] coordinates {(0,0.55) (-0.2,0.6) (-0.25,0.7) (-0.2,0.93) (-0.1, 1) (0,1.02) (0.1, 1) (0.2,0.93) (0.25,0.7) (0.2,0.6) (0,0.55)};
 \draw [<-,>=Stealth]  (0.195,0.95) -- (0.23,0.85);
\node[] at (-0.17, 0.75) {\fontsize{11}{11} $i$};
\node[] at (-0.25, -0.28) {\fontsize{11}{11} $-i$};
\filldraw (0,0.75) circle (0.75pt);
\filldraw (0,-0.3) circle (0.75pt);
\node[] at (0.3, 0.95) {\fontsize{9}{9} $\mathcal C$};
\draw [] (0.99, 1.14) -- (0.99, 1.26);
\draw [] (0.99, 1.14) -- (1.1, 1.14);
\node[] at (1.02, 1.213) {\fontsize{6}{6} $k$};
\end{tikzpicture}
\caption{The closed contour $\mathcal C$ for the integrals involving the boundary values.}
\label{contour-c-f}
\end{figure}
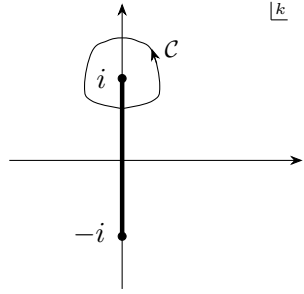

In particular, observe that, since $x>0$, the entire function $e^{ikx}$ is bounded for $\text{Im}(k)\geqslant 0$. Moreover, the functions $e^{i\omega(t-\tau)}+e^{-i\omega(t-\tau)}$ and $\frac{1}{\mu} \left[e^{i\omega(t-\tau)}+e^{-i\omega(t-\tau)}\right]$, which arise in the integrands of \eqref{serre-ir-U=0} through the relevant transforms of the boundary values,   are analytic and bounded for all $k\neq \pm i$. Hence, using Cauchy's theorem and Jordan's lemma from complex analysis (e.g. see Lemma 4.2.2 in \cite{af2003}), we are able to deform the contours of integration of the boundary values integrals in \eqref{serre-ir-U=0} from $\mathbb R$ to the closed contour $\mathcal C$ encircling $i$ (see Figure \ref{contour-c-f}).\footnote{This deformation is inspired by the one of \cite{vd2013} for the BBM equation; here, however, we have the additional complication of a dispersion relation with branching.} We emphasize that Jordan's lemma can be employed because of the  uniform (in $\arg(k)$) decay of the quantities $\frac{1}{1+k^2}$,  $\frac{1}{\mu(k)}$ and $\frac{k}{1+k^2}$ as $|k|\to \infty$. Hence, the integral representations \eqref{serre-ir-U=0} can be written as
\sss{\label{serre-ir-U=0-def}
\ddd{
r(x, t) 
&= \frac{1}{2\pi}  \int_{k\in \mathbb R} e^{ikx} \, \frac 12 \left[ \left(e^{i\omega t}+e^{-i\omega t}\right) \hat r_0(k) - \mu(k)  \left(e^{i\omega t}-e^{-i\omega t}\right) \hat q_0(k) \right]  dk
\nn\\
&\quad
+ \frac{1}{2\pi} \int_{k\in \mathcal C} e^{ikx} \, \frac{1}{2\mu(k)} \, \bigg\{ - \left[e^{i\omega t} \,  \widetilde h_0(\omega, t) - e^{-i\omega t} \, \widetilde h_0(-\omega, t)\right] 
-  \left[e^{i\omega t} \,  \widetilde{h_0''}(\omega, t) - e^{-i\omega t} \,  \widetilde{h_0''}(-\omega, t)\right]  
\nn\\
&\quad
+  \mu(k) \left[ e^{i\omega t} \,  \widetilde g_0(\omega, t) + e^{-i\omega t} \,  \widetilde g_0(-\omega, t) \right] 
+  ik   \left[e^{i\omega t} \,  \widetilde{g_0'}(\omega, t) - e^{-i\omega t} \,  \widetilde{g_0'}(-\omega, t)\right] \bigg\} \, dk
}
and
\ddd{
& q(x, t) 
= \frac{1}{2\pi} \int_{k\in \mathbb R} e^{ikx} \, \frac 12 \left[ \left(e^{i\omega t}+e^{-i\omega t}\right) \hat q_0(k)  
- \frac{1}{\mu(k)}  \left(e^{i\omega t}-e^{-i\omega t}\right) \hat r_0(k) \right] dk
\nn\\
&\quad
+ \frac{1}{2\pi} \int_{k\in \mathcal C} e^{ikx} \, \frac{1}{2\left(1+k^2\right)} \, \bigg\{  \left[ e^{i\omega t} \, \widetilde h_0(\omega, t) + e^{-i\omega t} \, \widetilde h_0(-\omega, t) \right]
+ \left[ e^{i\omega t} \,  \widetilde{h_0''}(\omega, t) + e^{-i\omega t} \,  \widetilde{h_0''}(-\omega, t)   \right]
\nn\\
&\quad
-  \mu(k) \left[ e^{i\omega t} \, \widetilde g_0(\omega, t) - e^{-i\omega t} \, \widetilde g_0(-\omega, t) \right]  
- ik \left[e^{i\omega t} \,  \widetilde{g_0'}(\omega, t) + e^{-i\omega t} \,  \widetilde{g_0'}(-\omega, t)\right] \bigg\}\, dk.
}
}

We are now ready to take advantage of the identities \eqref{serre-gr-U=0-sym}. Rearranging \eqref{serre-gr-U=0-sym-a}, we have
\sss{\label{lcbous-elim2}
\ddd{
&\quad
\frac{1}{2 \mu(k)} \, \bigg\{ \left[ e^{i\omega t} \,  \widetilde h_0(\omega, t) - e^{-i\omega t} \, \widetilde h_0(-\omega, t)\right]
+
 \left[e^{i\omega t} \,  \widetilde{h_0''}(\omega, t) - e^{-i\omega t} \,  \widetilde{h_0''}(-\omega, t)\right] \bigg\}
\\
&=
\hat r(-k, t) 
- \frac 12  \left[ \left(e^{i\omega t}+e^{-i\omega t}\right) \hat r_0(-k) 
+  \mu(k)  \left(e^{i\omega t}-e^{-i\omega t}\right) \hat q_0(-k) \right]
\nn\\
&\quad
- \frac{1}{2\mu(k)} \, \bigg\{ \mu(k) \left[ e^{i\omega t} \,  \widetilde g_0(\omega, t) + e^{-i\omega t} \,  \widetilde g_0(-\omega, t) \right]
+ ik  \left[e^{i\omega t} \,  \widetilde{g_0'}(\omega, t) - e^{-i\omega t} \,  \widetilde{g_0'}(-\omega, t)\right] \bigg\}.
\nn
}
Moreover, \eqref{serre-gr-U=0-sym-b} can be written as
\ddd{
&\quad
\frac{1}{2\left(1+k^2\right)} \, \bigg\{ \left[ e^{i\omega t} \, \widetilde h_0(\omega, t) + e^{-i\omega t} \, \widetilde h_0(-\omega, t)\right]
+ \left[e^{i\omega t} \,  \widetilde{h_0''}(\omega, t) + e^{-i\omega t} \,  \widetilde{h_0''}(-\omega, t)\right] \bigg\}
\\
&=
\hat q(-k, t) 
- \frac 12 \left[ \left(e^{i\omega t}+e^{-i\omega t}\right) \hat q_0(-k) + \frac{1}{\mu(k)}  \left(e^{i\omega t}-e^{-i\omega t}\right) \hat r_0(-k) \right]
\nn\\
&\quad
- \frac{1}{2\left(1+k^2\right)} \, \bigg\{ \mu(k) \left[ e^{i\omega t} \, \widetilde g_0(\omega, t) - e^{-i\omega t} \, \widetilde g_0(-\omega, t)\right]
+  ik  \left[e^{i\omega t} \,  \widetilde{g_0'}(\omega, t) + e^{-i\omega t} \,  \widetilde{g_0'}(-\omega, t)\right] \bigg\}.
\nn
}
}
Since the expressions \eqref{lcbous-elim2} are valid for all $\text{Im}(k) \geqslant 0$ with $k\neq i$, we employ them for $k\in \mathcal C$ and combine them with the integral representations \eqref{serre-ir-U=0-def} to obtain
\sss{\label{serre-ir-U=0-def-2}
\ddd{
r(x, t) 
&= \frac{1}{2\pi}  \int_{k\in \mathbb R} e^{ikx} \, \frac 12 \left[ \left(e^{i\omega t}+e^{-i\omega t}\right) \hat r_0(k) - \mu(k)  \left(e^{i\omega t}-e^{-i\omega t}\right) \hat q_0(k) \right]  dk
\nn\\
&\quad
+ \frac{1}{2\pi} \int_{k\in \mathcal C} e^{ikx}  \, \frac 12  \left[ \left(e^{i\omega t}+e^{-i\omega t}\right) \hat r_0(-k) 
+  \mu(k)  \left(e^{i\omega t}-e^{-i\omega t}\right) \hat q_0(-k) \right]   dk
\nn\\
&\quad
+ \frac{1}{2\pi} \int_{k\in \mathcal C} e^{ikx} \, \bigg\{ -\hat r(-k, t) +
  \left[ e^{i\omega t} \,  \widetilde g_0(\omega, t) + e^{-i\omega t} \,  \widetilde g_0(-\omega, t) \right] 
\nn\\
&\hskip 3.25cm
+ \frac{ik}{\mu(k)}    \left[e^{i\omega t} \,  \widetilde{g_0'}(\omega, t) - e^{-i\omega t} \,  \widetilde{g_0'}(-\omega, t)\right] \bigg\} \, dk
}
and
\ddd{
 q(x, t) 
&= \frac{1}{2\pi} \int_{k\in \mathbb R} e^{ikx} \, \frac 12 \left[ \left(e^{i\omega t}+e^{-i\omega t}\right) \hat q_0(k)  
- \frac{1}{\mu(k)}  \left(e^{i\omega t}-e^{-i\omega t}\right) \hat r_0(k) \right] dk
\nn\\
&\quad
- \frac{1}{2\pi} \int_{k\in \mathcal C} e^{ikx} \,  \frac 12 \left[ \left(e^{i\omega t}+e^{-i\omega t}\right) \hat q_0(-k) + \frac{1}{\mu(k)}  \left(e^{i\omega t}-e^{-i\omega t}\right) \hat r_0(-k) \right]   dk
\nn\\
&\quad
-\frac{1}{2\pi} \int_{k\in \mathcal C} e^{ikx} \,  \bigg\{ \hat q(-k, t)  + \frac{1}{\mu(k)} \left[ e^{i\omega t} \, \widetilde g_0(\omega, t) - e^{-i\omega t} \, \widetilde g_0(-\omega, t) \right] 
\nn\\
&\hskip 3.2cm 
+ \frac{ik}{1+k^2} \left[e^{i\omega t} \,  \widetilde{g_0'}(\omega, t) + e^{-i\omega t} \,  \widetilde{g_0'}(-\omega, t)\right] \bigg\}\, dk.
}
}

Of course, the above expressions still involve unknown quantities, namely the transforms $\hat r(-k, t)$ and $\hat q(-k, t)$. However, both of these transforms, as well as the exponential $e^{ikx}$, are analytic in the upper half-plane. Hence, by Cauchy's theorem inside the region enclosed by $\mathcal C$ we conclude that
\eee{
\int_{k\in \mathcal C} e^{ikx} \, \hat r(-k, t)  dk =   \int_{k\in \mathcal C} e^{ikx} \, \hat q(-k, t)  dk = 0
}
for all $x, t$. Therefore, we arrive at the following UTM solution formula for the linearized classical Boussinesq IBVP \eqref{lbous-ibvp}:
\sss{\label{lcbous-sol-fin}
\ddd{\label{lcbous-sol-fin-r}
r(x, t) 
&= \frac{1}{2\pi}  \int_{k\in \mathbb R} e^{ikx} \, \frac 12 \left[ \left(e^{i\omega t}+e^{-i\omega t}\right) \hat r_0(k) - \mu(k)  \left(e^{i\omega t}-e^{-i\omega t}\right) \hat q_0(k) \right]  dk
\nn\\
&\quad
+ \frac{1}{2\pi} \int_{k\in \mathcal C} e^{ikx}  \, \frac 12  \left[ \left(e^{i\omega t}+e^{-i\omega t}\right) \hat r_0(-k) 
+  \mu(k)  \left(e^{i\omega t}-e^{-i\omega t}\right) \hat q_0(-k) \right]   dk
\nn\\
&\quad
+ \frac{1}{2\pi} \int_{k\in \mathcal C} e^{ikx} \, \bigg\{ 
  \left[ e^{i\omega t} \,  \widetilde g_0(\omega, t) + e^{-i\omega t} \,  \widetilde g_0(-\omega, t) \right] 
\nn\\
&\hskip 3.25cm
+ i  \omega(k)   \left[e^{i\omega t} \,  \widetilde{g_0'}(\omega, t) - e^{-i\omega t} \,  \widetilde{g_0'}(-\omega, t)\right] \bigg\} \, dk
}
and
\ddd{\label{lcbous-sol-fin-q}
 q(x, t) 
&= \frac{1}{2\pi} \int_{k\in \mathbb R} e^{ikx} \, \frac 12 \left[ \left(e^{i\omega t}+e^{-i\omega t}\right) \hat q_0(k)  
- \frac{1}{\mu(k)}  \left(e^{i\omega t}-e^{-i\omega t}\right) \hat r_0(k) \right] dk
\nn\\
&\quad
- \frac{1}{2\pi} \int_{k\in \mathcal C} e^{ikx} \,  \frac 12 \left[ \left(e^{i\omega t}+e^{-i\omega t}\right) \hat q_0(-k) + \frac{1}{\mu(k)}  \left(e^{i\omega t}-e^{-i\omega t}\right) \hat r_0(-k) \right]   dk
\nn\\
&\quad
-\frac{1}{2\pi} \int_{k\in \mathcal C} e^{ikx} \,  \bigg\{  \frac{1}{\mu(k)} \left[ e^{i\omega t} \, \widetilde g_0(\omega, t) - e^{-i\omega t} \, \widetilde g_0(-\omega, t) \right] 
\nn\\
&\hskip 3.2cm 
+ \frac{ik}{1+k^2} \left[e^{i\omega t} \,  \widetilde{g_0'}(\omega, t) + e^{-i\omega t} \,  \widetilde{g_0'}(-\omega, t)\right] \bigg\}\, dk.
}
}
We emphasize that, as explained in Remark \ref{bcut-r}, the integrands in the above formula are \textit{analytic} functions of $k$ in the respective domains of integration despite the fact that they involve the branched function $\omega(k)$.

\begin{remark}[Deformation back to $\mathbb R$]
Thanks to analyticity and exponential decay, it it is possible to deform the contour of integration $\mathcal C$ in formulae \eqref{lcbous-sol-fin} back to $\mathbb R$. However, the resulting expression is \textit{not} uniformly convergent at the boundary $x=0$, and hence it is \textit{not} suitable for explicitly verifying that the UTM formulae \eqref{lcbous-sol-fin} indeed satisfy IBVP \eqref{lbous-ibvp}. This fact is clearly illustrated by the computations of Section \ref{ver-s}. 
\end{remark}

\begin{remark}[Other types of boundary conditions]
The elimination procedure performed in this section works in the same way for other types of admissible boundary data. For example, instead of the Dirichlet condition \eqref{lbous-bc}, one could prescribe the Neumann datum $q_x(0, t) = g_1(t)$ which, by integrating \eqref{g1-h0-rel} and employing the compatibility condition $h_0(0) = r_0(0)$, is  equivalent to the prescription of the Dirichlet datum  $r(0, t) = h_0(t)$.  
\end{remark}

\begin{remark}[Solution for the forced linear problem]
Thanks to the Duhamel principle, the UTM solution formula \eqref{lcbous-sol-fin} can be easily adapted to the case of the forced counterpart of the linearized classical Boussinesq IBVP \eqref{lbous-ibvp}. The resulting formula can then be employed for studying the well-posedness of the nonlinear classical Boussinesq system \eqref{cbous} on the half-line via contraction mapping techniques.
\end{remark}

\section{Revisiting the problem as a single equation}
\label{one-eq-s}

In the previous section, we solved IBVP \eqref{lbous-ibvp} for the linearized classical Boussinesq system using the recently introduced extension of UTM to systems \cite{dgsv2018}. This extension provides a general method that can in principle be applied  to any linear system of evolution equations. However, specifically in the case of IBVP \eqref{lbous-ibvp}, it is also possible to convert the problem into one involving a single equation and hence solve it via the standard UTM. Indeed, differentiating the first equation of system \eqref{lbous-sys} with respect to $x$ and the second one with respect to $t$, we have\footnote{Recall that we are working with smooth functions and hence we are allowed to interchange the order of partial derivatives.}
\eee{\label{r-q-eq}
r_{xt} +  q_{xx} = 0, \quad q_{tt}+   r_{xt} - q_{xxtt} = 0.
}
Thus, we can eliminate $r$ and obtain a single equation for $q$:
\begin{equation}\label{q-eq}
q_{tt} - q_{xx} - q_{xxtt} = 0.
\end{equation}
In the remaining of this section, we will employ UTM to solve equation \eqref{q-eq} in terms of the initial and boundary data of problem \eqref{lbous-ibvp}. Once an explicit solution formula for $q$ is obtained, it will be straightforward to deduce a corresponding formula for $r$ since a simple integration of the first of equations \eqref{r-q-eq} yields
\eee{
r(x, t) = r_0(x) -    \int_{\tau=0}^t q_x(x, \tau) d\tau
}
with $r_0(x)$ being the initial datum for $r$ prescribed in problem \eqref{lbous-ibvp}.

Applying the half-line Fourier transform \eqref{ft-def} to equation \eqref{q-eq}, we find
\eee{\label{r-ode}
\p_t^2 \hat q(k, t) + \frac{k^2}{1+k^2} \, \hat q(k, t)
=
-\frac{1}{1+k^2}   \left[g_1(t) + ik g_0(t) + g_1''(t) + ik g_0''(t)\right],
}
where  $\text{Im}(k) \leqslant 0$ with $k\neq -i$ and we have denoted, as usual, $g_0(t) = q(0, t)$ and $g_1(t) = q_x(0, t)$. 
The second-order ordinary differential equation \eqref{r-ode} can be solved via variation of parameters.

In particular,   the general solution to the homogeneous counterpart of \eqref{r-ode} is
\eee{\label{rhom}
\hat q_h(k, t) = c_1(k) e^{i\omega t} + c_2(k) e^{-i\omega t}
}
with $\omega$ defined by \eqref{om-def-eq}.
The constants (with respect to $t$) $c_1(k)$ and $c_2(k)$ in the homogeneous solution \eqref{rhom} can be computed by enforcing the initial conditions of problem \eqref{lbous-ibvp}. First, note that the initial condition $q(x, 0) = q_0(x)$  readily implies 
\eee{\label{rh-ic1}
\hat q(k, 0) = \hat q_0(k).
}
Furthermore,  the second of equations \eqref{lbous-sys} evaluated at $t=0$ and combined with the initial condition $r(x, 0) =r_0(x)$ yields
\eee{
q_t(x, 0) -  q_{xxt}(x, 0) = -  r_0'(x)
}
and, therefore, taking the half-line Fourier transform \eqref{ft-def} we obtain
\eee{\label{rh-ic2}
\hat q_t(k, 0) = \frac{1}{1+k^2} \left[ r_0(0) -ik  \hat r_0(k)\ - g_1'(0) - ik g_0'(0)\right] =: \hat q_1(k).
}
The conditions \eqref{rh-ic1} and \eqref{rh-ic2} must be satisfied by the homogeneous solution formula \eqref{rhom}. For this, we must have
\eee{\label{c1c2-def}
c_1(k) = \frac 12 \left[ \hat q_0(k) + \frac{\hat q_1(k)}{i\omega} \right], \quad  c_2(k) = \frac 12 \left[\hat q_0(k) - \frac{\hat q_1(k)}{i\omega} \right].
}

Furthermore,  variation of parameters yields a particular solution of \eqref{r-ode} in the form
\eee{
\hat q_p(k, t) 
=
e^{i\omega t} \int_{\tau=0}^t \frac{e^{-i\omega \tau} \varphi(\tau)}{2i\omega} \, d\tau
-
e^{-i\omega t} \int_{\tau=0}^t \frac{e^{i\omega \tau} \varphi(\tau)}{2i\omega} \, d\tau,
}
where $\varphi(t) $ is simply the forcing on the right-hand side of \eqref{r-ode}:
\eee{
\varphi(t) = \varphi(k, t) = -\frac{1}{1+k^2} \left[g_1(t) + ik g_0(t) + g_1''(t) + ik g_0''(t)\right].
}

Therefore, recalling the  notation $\widetilde f(\omega, t) = \int_{\tau=0}^t e^{-i\omega \tau} f(\tau) d\tau$ and noting that $\hat q = \hat q_h + \hat q_p$, we overall find the solution to \eqref{r-ode} as
\ddd{\label{r-eq-gr}
\hat q(k, t) 
&=
\frac 12 \left(e^{i\omega t} +e^{-i\omega t}\right) \hat q_0(k) 
+
\frac{1}{2i\omega} \left(e^{i\omega t} - e^{-i\omega t}\right) \hat q_1(k) 
\nn\\
&\quad
+\frac{1}{2i\omega \left(1+k^2\right)} \,
\bigg\{
-e^{i\omega t} \left[
 \widetilde g_1(\omega, t) + \widetilde{g_1''}(\omega, t)  + ik \widetilde g_0(\omega, t) + ik \widetilde{g_0''}(\omega, t)  \right]
\nn\\
&\quad
+ e^{-i\omega t} 
\left[
  \widetilde g_1(-\omega, t) + \widetilde{g_1''}(-\omega, t) + ik \widetilde g_0(-\omega, t) + ik \widetilde{g_0''}(-\omega, t) 
\right]
\bigg\},
} 
valid for $\text{Im}(k) \leqslant 0$ with $k\neq -i$.
In fact, integrating by parts we have
$$
\widetilde{g_0''}(\pm \omega, t) 
=
e^{\mp i\omega t} g_0'(t) - g_0'(0) \pm i\omega \widetilde{g_0'}(\pm \omega, t) 
$$
and similarly for $\widetilde{g_1''}(\pm \omega, t)$.
Thus, substituting also for $\hat q_1$ via \eqref{rh-ic2}, we can write \eqref{r-eq-gr} in the form
\ddd{\label{r-eq-gr-2}
\hat q(k, t) 
&=
\frac 12 \left(e^{i\omega t} +e^{-i\omega t}\right) \hat q_0(k) 
+
\frac{1}{2i\omega \left(1+k^2\right)} \left(e^{i\omega t} - e^{-i\omega t}\right)   \left[ r_0(0) -ik  \hat r_0(k)\right] 
\\
&\quad
-\frac{1}{2i\omega \left(1+k^2\right)} \, 
\bigg\{
 \left[e^{i\omega t} \, \widetilde g_1(\omega, t)  - e^{-i\omega t} \, \widetilde g_1(-\omega, t)\right] 
+
i\omega \left[e^{i\omega t} \, \widetilde{g_1'}(\omega, t)  + e^{-i\omega t} \, \widetilde{g_1'}(-\omega, t)\right]
\nn\\
&\hskip 2.8cm
+ ik  \left[e^{i\omega t} \, \widetilde g_0(\omega, t)  - e^{-i\omega t} \, \widetilde g_0(-\omega, t)\right] 
- k\omega \left[e^{i\omega t} \, \widetilde{g_0'}(\omega, t)  + e^{-i\omega t} \, \widetilde{g_0'}(-\omega, t)\right]
\bigg\}.
\nn
} 

Expression \eqref{r-eq-gr-2} is the global relation corresponding to equation \eqref{q-eq}. Since it is valid for $\text{Im}(k) \leqslant 0$ with $k\neq -i$, it can be combined with the Fourier inversion formula \eqref{ift-def} to yield an integral representation for $q$. This integral representation will involve two boundary values, $g_0$ and  $g_1$. Since only the former one is prescribed as a boundary condition in problem \eqref{lbous-ibvp},  the latter one will have to be eliminated. This is achieved by using the transformation $k\mapsto -k$, which is a symmetry of $\omega$, and then by exploiting analyticity and Cauchy's theorem. The whole procedure is just like the one presented in detail in Section \ref{sol-s} and so we do not repeat it here. In fact, it is easy to convert  the global relation \eqref{r-eq-gr-2} to the global relation \eqref{serre-gr-U=0-r} that we obtained earlier in Section \ref{sol-s} using the UTM for systems approach. 
 
 Indeed,  integrating by parts and recalling \eqref{g1-h0-rel}, we find
\ddd{
&\quad
e^{i\omega t} \, \widetilde{g_1}(\omega, t) - e^{-i\omega t} \, \widetilde{g_1}(-\omega, t)
=
-e^{i\omega t} \, \widetilde{h_0'}(\omega, t) + e^{-i\omega t} \, \widetilde{h_0'}(-\omega, t)
\nn\\
&=
-e^{i\omega t} \left[e^{-i\omega t} \, h_0(t) - h_0(0) + i\omega \widetilde h_0(\omega, t) \right]
+e^{-i\omega t} \left[e^{i\omega t} h_0(t) - h_0(0) - i\omega \widetilde h_0(-\omega, t) \right]
\nn\\
&=
\left(e^{i\omega t} - e^{-i\omega t}\right) r_0(0)
-i\omega \left[ e^{i\omega t} \, \widetilde h_0(\omega, t) + e^{-i\omega t} \, \widetilde h_0(-\omega, t) \right],
}
where the last equality follows from the compatibility at the origin: $h_0(0) := r(0, 0) =:  r_0(0)$.
In addition, by \eqref{g1-h0-rel} we have
\eee{
e^{i\omega t} \, \widetilde{g_1'}(\omega, t) + e^{-i\omega t} \, \widetilde{g_1'}(-\omega, t)
=
- e^{i\omega t} \, \widetilde{h_0''}(\omega, t) - e^{-i\omega t} \, \widetilde{h_0''}(-\omega, t).
}
Therefore, \eqref{r-eq-gr-2} can also be written in the form
%
\ddd{\label{r-eq-gr-3}
\hat q(k, t) 
&=
\frac 12 \left(e^{i\omega t} +e^{-i\omega t}\right) \hat q_0(k) 
-
\frac{1}{2\mu(k)} \left(e^{i\omega t} - e^{-i\omega t}\right)   \hat r_0(k) 
\nn\\
&\quad
+\frac{1}{2\left(1+k^2\right)} \,
\bigg\{
\left[ e^{i\omega t} \, \widetilde h_0(\omega, t) + e^{-i\omega t} \, \widetilde h_0(-\omega, t) \right]
+
 \left[e^{i\omega t} \, \widetilde{h_0''}(\omega, t)  + e^{-i\omega t} \, \widetilde{h_0''}(-\omega, t)\right]
\nn\\
&\quad
-
 \mu(k)  \left[e^{i\omega t} \, \widetilde g_0(\omega, t)  - e^{-i\omega t} \, \widetilde g_0(-\omega, t)\right] 
-
ik \left[e^{i\omega t} \, \widetilde{g_0'}(\omega, t)  + e^{-i\omega t} \, \widetilde{g_0'}(-\omega, t)\right]
\bigg\},
}
which is precisely the global relation \eqref{serre-gr-U=0-r} that we obtained in Section \ref{sol-s} via the UTM for systems approach. Therefore, from this point onwards, following the elimination procedure of Section \ref{sol-s} we arrive once again at formula \eqref{lcbous-sol-fin-q} for $q$. Once again, we emphasize that it is possible to obtain formula \eqref{lcbous-sol-fin-q} directly from the global relation \eqref{r-eq-gr-2}, i.e. without converting that global relation into \eqref{serre-gr-U=0-r}.

\section{Explicit verification of the novel solution formula}
\label{ver-s}

Since the derivations of the previous sections were performed under the assumption of existence of solution, we shall now verify that the resulting formulae do indeed satisfy   IBVP \eqref{lbous-ibvp}. 

We begin with the  linearized classical Boussinesq system \eqref{lbous-sys}. Differentiating formula \eqref{lcbous-sol-fin-r} with respect to $t$, we have
\ddd{\label{lcbous-sol-fin-r-t}
r_t(x, t) 
&= \frac{1}{2\pi}  \int_{k\in \mathbb R} e^{ikx} \, \frac{i\omega}2 \left[ \left(e^{i\omega t}-e^{-i\omega t}\right) \hat r_0(k) - \mu(k)  \left(e^{i\omega t}+e^{-i\omega t}\right) \hat q_0(k) \right]  dk
\nn\\
&\quad
+ \frac{1}{2\pi} \int_{k\in \mathcal C} e^{ikx}  \, \frac{i\omega}2  \left[ \left(e^{i\omega t}-e^{-i\omega t}\right) \hat r_0(-k) 
+  \mu(k)  \left(e^{i\omega t}+e^{-i\omega t}\right) \hat q_0(-k) \right]   dk
\nn\\
&\quad
+ \frac{1}{2\pi} \int_{k\in \mathcal C} e^{ikx} \, i\omega \, \bigg\{ 
  \left[ e^{i\omega t} \,  \widetilde g_0(\omega, t) - e^{-i\omega t} \,  \widetilde g_0(-\omega, t) \right] 
\nn\\
&\hskip 3.65cm
+ i  \omega(k)   \left[e^{i\omega t} \,  \widetilde{g_0'}(\omega, t) + e^{-i\omega t} \,  \widetilde{g_0'}(-\omega, t)\right] \bigg\} \, dk
\nn\\
&\quad
+ \frac{1}{2\pi} \int_{k\in \mathcal C} e^{ikx} \, \bigg\{ 
  \left[  g_0(t) +  g_0(t) \right] 
+ i \omega(k)   \left[g_0'(t) - g_0'(t)\right] \bigg\} \, dk,
}
while taking the derivative of \eqref{lcbous-sol-fin-q} with respect to $x$ gives
\ddd{\label{lcbous-sol-fin-q-x}
 q_x(x, t) 
&= \frac{1}{2\pi} \int_{k\in \mathbb R} e^{ikx} \, \frac{ik}2 \left[ \left(e^{i\omega t}+e^{-i\omega t}\right) \hat q_0(k)  
- \frac{1}{\mu(k)}  \left(e^{i\omega t}-e^{-i\omega t}\right) \hat r_0(k) \right] dk
\nn\\
&\quad
- \frac{1}{2\pi} \int_{k\in \mathcal C} e^{ikx} \,  \frac{ik}2 \left[ \left(e^{i\omega t}+e^{-i\omega t}\right) \hat q_0(-k) + \frac{1}{\mu(k)}  \left(e^{i\omega t}-e^{-i\omega t}\right) \hat r_0(-k) \right]   dk
\nn\\
&\quad
-\frac{1}{2\pi} \int_{k\in \mathcal C} e^{ikx} \, ik \,  \bigg\{  \frac{1}{\mu(k)} \left[ e^{i\omega t} \, \widetilde g_0(\omega, t) - e^{-i\omega t} \, \widetilde g_0(-\omega, t) \right] 
\nn\\
&\hskip 3.6cm 
+ \frac{ik}{1+k^2} \left[e^{i\omega t} \,  \widetilde{g_0'}(\omega, t) + e^{-i\omega t} \,  \widetilde{g_0'}(-\omega, t)\right] \bigg\}\, dk.
}
Note, importantly, that the last integral in \eqref{lcbous-sol-fin-r-t} is zero, since the function $e^{ikx}$ is analytic inside the region enclosed by $\mathcal C$  and the part of the integrand involving $\omega(k)$ vanishes identically. Thus, recalling that $\omega \mu = k$, we see that $r_t + q_x = 0$, i.e. the first equation of system \eqref{lbous-sys} is satisfied.

Furthermore, using Cauchy's residue theorem we compute
$$
\int_{k\in \mathcal C} e^{ikx} \,   \frac{ik}{1+k^2} \,  dk = - \pi  e^{-x}
$$
and hence
\ddd{\label{lcbous-sol-fin-qt}
q_t(x, t) 
&= \frac{1}{2\pi} \int_{k\in \mathbb R} e^{ikx} \, \frac{i\omega}2 \left[ \left(e^{i\omega t}-e^{-i\omega t}\right) \hat q_0(k)  
- \frac{1}{\mu(k)}  \left(e^{i\omega t}+e^{-i\omega t}\right) \hat r_0(k) \right] dk
\nn\\
&\quad
- \frac{1}{2\pi} \int_{k\in \mathcal C} e^{ikx} \,  \frac{i\omega}2 \left[ \left(e^{i\omega t}-e^{-i\omega t}\right) \hat q_0(-k) + \frac{1}{\mu(k)}  \left(e^{i\omega t}+e^{-i\omega t}\right) \hat r_0(-k) \right]   dk
\nn\\
&\quad
-\frac{1}{2\pi} \int_{k\in \mathcal C} e^{ikx} \, i\omega\,  \bigg\{  \frac{1}{\mu(k)} \left[ e^{i\omega t} \, \widetilde g_0(\omega, t) + e^{-i\omega t} \, \widetilde g_0(-\omega, t) \right] 
\nn\\
&\hskip 3.6cm 
+ \frac{ik}{1+k^2} \left[e^{i\omega t} \,  \widetilde{g_0'}(\omega, t) - e^{-i\omega t} \,  \widetilde{g_0'}(-\omega, t)\right] \bigg\}\, dk
+e^{-x} g_0'(t)
}
and 
\ddd{\label{lcbous-sol-fin-qxxt}
q_{xxt}(x, t) 
&= \frac{1}{2\pi} \int_{k\in \mathbb R} e^{ikx} \, \frac{-ik^2\omega}2 \left[ \left(e^{i\omega t}-e^{-i\omega t}\right) \hat q_0(k)  
- \frac{1}{\mu(k)}  \left(e^{i\omega t}+e^{-i\omega t}\right) \hat r_0(k) \right] dk
\\
&\quad
- \frac{1}{2\pi} \int_{k\in \mathcal C} e^{ikx} \,  \frac{-ik^2\omega}2 \left[ \left(e^{i\omega t}-e^{-i\omega t}\right) \hat q_0(-k) + \frac{1}{\mu(k)}  \left(e^{i\omega t}+e^{-i\omega t}\right) \hat r_0(-k) \right]   dk
\nn\\
&\quad
-\frac{1}{2\pi} \int_{k\in \mathcal C} e^{ikx} \left(-ik^2\omega\right)  \bigg\{  \frac{1}{\mu(k)} \left[ e^{i\omega t} \, \widetilde g_0(\omega, t) + e^{-i\omega t} \, \widetilde g_0(-\omega, t) \right] 
\nn\\
&\hskip 4.7cm 
+ \frac{ik}{1+k^2} \left[e^{i\omega t} \,  \widetilde{g_0'}(\omega, t) - e^{-i\omega t} \,  \widetilde{g_0'}(-\omega, t)\right] \bigg\}\, dk
+e^{-x} g_0'(t).
\nn
}
Moreover,
\ddd{\label{lcbous-sol-fin-r-x}
r_x(x, t) 
&= \frac{1}{2\pi}  \int_{k\in \mathbb R} e^{ikx} \, \frac{ik}2 \left[ \left(e^{i\omega t}+e^{-i\omega t}\right) \hat r_0(k) - \mu(k)  \left(e^{i\omega t}-e^{-i\omega t}\right) \hat q_0(k) \right]  dk
\nn\\
&\quad
+ \frac{1}{2\pi} \int_{k\in \mathcal C} e^{ikx}  \, \frac{ik}2  \left[ \left(e^{i\omega t}+e^{-i\omega t}\right) \hat r_0(-k) 
+  \mu(k)  \left(e^{i\omega t}-e^{-i\omega t}\right) \hat q_0(-k) \right]   dk
\nn\\
&\quad
+ \frac{1}{2\pi} \int_{k\in \mathcal C} e^{ikx} \, ik \, \bigg\{ 
  \left[ e^{i\omega t} \,  \widetilde g_0(\omega, t) + e^{-i\omega t} \,  \widetilde g_0(-\omega, t) \right] 
\nn\\
&\hskip 3.65cm
+ \frac{ik}{\mu(k)}    \left[e^{i\omega t} \,  \widetilde{g_0'}(\omega, t) - e^{-i\omega t} \,  \widetilde{g_0'}(-\omega, t)\right] \bigg\} \, dk.
}
Therefore, observing that $\omega - k\mu + k^2 \omega = 0$, we deduce that $q_t + r_x - q_{xxt} =0$, as required by the second component of system \eqref{lbous-sys}.

Next, we verify the initial conditions \eqref{lbous-ic}. Evaluating \eqref{lcbous-sol-fin} at $t=0$, we have
\eee{\label{lcbous-sol-fin-r-t=0}
r(x, 0) 
= \frac{1}{2\pi}  \int_{k\in \mathbb R} e^{ikx} \,   \hat r_0(k)   dk
+ \frac{1}{2\pi} \int_{k\in \mathcal C} e^{ikx}  \,   \hat r_0(-k)   dk
}
and
\ddd{\label{lcbous-sol-fin-q-t=0}
 q(x, 0) 
&= \frac{1}{2\pi} \int_{k\in \mathbb R} e^{ikx} \,   \hat q_0(k)    dk
- \frac{1}{2\pi} \int_{k\in \mathcal C} e^{ikx} \,    \hat q_0(-k)   dk.
}
In view of the Fourier inversion formula \eqref{ift-def}, the first integrals in the above expressions are simply $r_0(x)$ and $q_0(x)$, respectively. Moreover, the second integrals are zero because of analyticity of the integrands inside the region enclosed by $\mathcal C$. Hence, the initial conditions are satisfied.

Finally, we verify the boundary condition \eqref{lbous-bc}. This task is a bit more challenging that the previous two verifications. We begin by observing that
$$
\hat q_0(k) := \int_{y=0}^\infty e^{-iky} q_0(y) dy = \frac{1}{ik} \left[q_0(0) + \what{q_0'}(k) \right]
$$
and, similarly,
$
\widetilde{g_0'}(\omega, t) 
=
e^{-i\omega t} g_0(t) - g_0(0) + i\omega \widetilde g_0(\omega, t).
$
Thus, noting that
$$
\int_{k\in \mathbb R} e^{ikx}   \left(e^{i\omega t} +e^{-i\omega t}\right) \hat q_0(k) dk
=
\dashint_{k\in \mathbb R} e^{ikx}   \left(e^{i\omega t} +e^{-i\omega t}\right) \hat q_0(k) dk
$$
 since the integral on the left-hand side exists without the need for taking the principal value $\dashint$, we have
\ddd{\label{r-eq-sol-ver1}
q(x, t)
&=
\frac{1}{2i\pi} \dashint_{k\in \mathbb R} e^{ikx}  \, \frac{\cos(\omega t)}{k}  \left[ q_0(0) + \what{q_0'}(k) \right]
 dk
+
\frac{1}{2i\pi} \int_{k\in \mathcal C} e^{ikx} 
\, \frac{\cos(\omega t)}{k}  \left[ q_0(0) + \what{q_0'}(-k) \right] dk
\nn\\
&\quad
+\frac{1}{2i\pi} \int_{k\in \mathbb R} e^{ikx} 
\, \frac{\sin(\omega t)}{\mu(k)} \,  \hat r_0(k)  dk
+
\frac{1}{2i\pi} \int_{k\in \mathcal C} e^{ikx} \, 
\frac{\sin(\omega t)}{\mu(k)}  \,  \hat r_0(-k) dk
\nn\\
&\quad
-\frac{1}{2\pi} \int_{k\in \mathcal C} e^{ikx} 
\frac{1}{\left(1+k^2\right)^{\frac 32}}
\left[e^{i\omega t} \, \widetilde g_0(\omega, t)  - e^{-i\omega t} \, \widetilde g_0(-\omega, t)\right]  dk
\nn\\
&\quad
+
\frac{g_0(t)}{i\pi} \int_{k\in \mathcal C} e^{ikx} \, \frac{k}{1 + k^2}\,  dk
-
\frac{g_0(0)}{2i\pi} \int_{k\in \mathcal C} e^{ikx} 
\frac{k}{1+k^2}
 \left(e^{i\omega t} + e^{-i\omega t}\right) dk.
}

Now, note that all the integrals in \eqref{r-eq-sol-ver1} whose contour is $\mathcal C$ are uniformly convergent and hence we can pass the limit $x\to 0$ inside them. The same is true for the integral along $\mathbb R$  involving $\sin(\omega t)$, since the integrand is of $O(k^{-2})$ as $|k|\to \infty$.\footnote{Indeed, we have $\omega \simeq 1$ as $|k|\to \infty$ and so $\sin(\omega t)$ is bounded. Moreover, $\mu(k) \simeq k$ and, finally, integrating by parts yields $\hat r_0(k) = \frac{1}{ik}r_0(0) + \frac{1}{ik} \what{r_0'}(k)$.} The first integral, i.e. the one along $\mathbb R$   involving $\cos(\omega t)$, will be discussed separately below. Furthermore, 
\ddd{
\int_{k\in \mathcal C}  
\frac{\sin(\omega t)}{\mu(k)}  \,  \hat r_0(-k) dk
&:=
\int_{y=0}^\infty r_0(y) \int_{k\in \mathcal C}  e^{iky} \, 
\frac{\sin(\omega t)}{\mu(k)}  \,  dk dy
\nn\\
&=
\int_{y=0}^\infty r_0(y) \int_{k\in \mathbb R}  e^{iky} \, 
\frac{\sin(\omega t)}{\mu(k)}  \,  dk dy
=
\int_{k\in \mathbb R}  
\frac{\sin(\omega t)}{\mu(k)}  \,  \hat r_0(-k) dk
}
by applying Cauchy's theorem and Jordan's lemma along the upper semicircle  of infinite radius in the complex $k$-plane.\footnote{Crucial for the application of Jordan's lemma is the uniform decay of the non-exponential part of the integrand thanks to $\frac{1}{\mu(k)}$.} 
Therefore,   \eqref{r-eq-sol-ver1} becomes
\ddd{
q(0, t)
&=
\frac{1}{2i\pi} \lim_{x\to 0} \dashint_{k\in \mathbb R} e^{ikx}  \, \frac{\cos(\omega t)}{k}  \left[ q_0(0) + \what{q_0'}(k) \right]
 dk
+
\frac{1}{2i\pi}   \int_{k\in \mathcal C}  \frac{\cos(\omega t)}{k}  \left[ q_0(0) + \what{q_0'}(-k) \right] dk
\nn\\
&\quad
+\frac{1}{2i\pi} \int_{k\in \mathbb R} 
\, \frac{\sin(\omega t)}{\mu(k)} \,  \hat r_0(k)  dk
+
\frac{1}{2i\pi} \int_{k\in \mathbb R}  
\frac{\sin(\omega t)}{\mu(k)}  \,  \hat r_0(-k) dk
\nn\\
&\quad
-\frac{1}{2\pi} \int_{k\in \mathcal C} 
\frac{1}{\left(1+k^2\right)^{\frac 32}}
\left[e^{i\omega t}\widetilde g_0(\omega, t)  - e^{-i\omega t} \widetilde g_0(-\omega, t)\right]  dk
\nn\\
&\quad
+
\frac{g_0(t)}{i\pi} \int_{k\in \mathcal C} 
\frac{k}{1+k^2}\,  dk
-
\frac{g_0(0)}{2i\pi} \int_{k\in \mathcal C} 
\frac{k}{1+k^2}
 \left(e^{i\omega t} + e^{-i\omega t}\right) dk,
}
and making the change of variable $k\mapsto -k$ to see that the two integrals involving $\hat r_0$ cancel (recall that $\mu(-k) = - \mu(k)$), we obtain
\ddd{\label{r-eq-sol-ver2}
q(0, t)
&=
\frac{1}{2i\pi} \lim_{x\to 0} \dashint_{k\in \mathbb R} e^{ikx}  \, \frac{\cos(\omega t)}{k}  \left[ q_0(0) + \what{q_0'}(k) \right]
 dk
+
\frac{1}{2i\pi}   \int_{k\in \mathcal C}  \frac{\cos(\omega t)}{k}  \left[ q_0(0) + \what{q_0'}(-k) \right] dk
\nn\\
&\quad
+ \frac{1}{i\pi} \int_{\tau=0}^t  g_0(\tau) \int_{k\in \mathcal C} 
\frac{\sin(\omega(t-\tau))}{\left(1+k^2\right)^{\frac 32}} \, dk d\tau
+
\frac{g_0(t)}{i\pi} \int_{k\in \mathcal C} 
\frac{k}{1+k^2}\,  dk
-
\frac{g_0(0)}{i\pi} \int_{k\in \mathcal C} 
\frac{k \cos(\omega t)}{1+k^2}  dk.
}

Next, we compute several integrals by exploiting the uniform convergence of Taylor series for $\cos$ and $\sin$ and using Cauchy's residue theorem.
First,  we have
\eee{\label{complex-int}
\int_{k\in \mathcal C}  \frac{\cos(\omega t)}{k}  \, q_0(0) dk
=
q_0(0) \sum_{j=0}^\infty (-1)^j \frac{t^{2j}}{(2j)!} \underbrace{\int_{k\in \mathcal C}  \frac{k^{2j-1}}{\left(1+k^2\right)^j} \,  dk}_{= \, i\pi \text{ for } j\in \mathbb N \text{ and } 0 \text{ for } j=0}
=
i\pi q_0(0) \left(\cos t -1\right),
}
where the above integral has been evaluated by using the standard complex analysis formula for the residue of a pole together with the Leibniz rule for the derivative of a product. Similarly,
\ddd{
\int_{k\in \mathcal C} 
\frac{\sin(\omega(t-\tau))}{\left(1+k^2\right)^{\frac 32}} \, dk 
&=
\sum_{j=0}^\infty (-1)^j \frac{\left(t-\tau\right)^{2j+1}}{\left(2j+1\right)!}
\underbrace{\int_{k\in \mathcal C} \frac{k^{2j+1}}{\left(1+k^2\right)^{j+2}} \, dk}_{= \, 0\, \forall j \in \mathbb N \cup \{0\}}
= 0.
}
Also, $\int_{k\in \mathcal C} 
\frac{k}{1+k^2}\,  dk = i\pi$ and, finally, 
\eee{
\int_{k\in \mathcal C} 
\frac{k \cos(\omega t)}{1+k^2}  dk
=
\sum_{j=0}^\infty (-1)^j \frac{t^{2j}}{(2j)!} \underbrace{\int_{k\in \mathcal C} 
\frac{k^{2j+1}}{\left(1+k^2\right)^{j+1}}  dk}_{= \, i\pi \, \forall j\in \mathbb N \cup \{0\}}
  = i\pi \cos t .
}
Substituting the above computations in \eqref{r-eq-sol-ver2} and recalling the compatibility condition $q_0(0) = g_0(0)$, we obtain
\ddd{\label{r-eq-sol-ver3}
q(0, t)
&=
\frac{1}{2i\pi} \lim_{x\to 0} \dashint_{k\in \mathbb R} e^{ikx}  \, \frac{\cos(\omega t)}{k}  \left[ g_0(0) + \what{q_0'}(k) \right] dk
+
\frac{1}{2i\pi}   \int_{k\in \mathcal C}  \frac{\cos(\omega t)}{k}  \, \what{q_0'}(-k)   dk
\nn\\
&\quad
-
\frac 12 g_0(0) \left(\cos t+1\right)
+
g_0(t).
}

Next, we discuss the first principal value integral in \eqref{r-eq-sol-ver3}. We have
\eee{
 \dashint_{k\in \mathbb R} e^{ikx}  \, \frac{\cos(\omega t)}{k}  \,   dk
=
\int_{k\in \Gamma} e^{ikx}  \, \frac{\cos(\omega t)}{k}  \,   dk
-
  \int_{k\in C_\varepsilon} e^{ikx}  \, \frac{\cos(\omega t)}{k}  \,   dk,
}
where $\Gamma$ is the closed, anti-clockwise contour consisting of $(-\infty, -\varepsilon]$, $C_\varepsilon$, $[\varepsilon, \infty)$ and $C_R$, with $C_\varepsilon$ being the upper semicircle of radius $\varepsilon$ centered at the origin and oriented clockwise, and with $C_R$ being the upper semicircle of radius $R\to \infty$ centered at the origin and oriented anti-clockwise. We note that the deformation from $(-\infty, -\varepsilon] \cup C_\varepsilon \cup [\varepsilon, \infty)$ to $\Gamma$ is possible due to the fact that the integral along $C_R$ vanishes thanks to Jordan's lemma.\footnote{Note here the importance of  the uniform decay of $\frac 1k$ in the integrand (recall that $\omega \simeq 1$ as $|k|\to \infty$ so $\cos(\omega t)$ is bounded at infinity).} Now, by Cauchy's theorem and \eqref{complex-int} we have
\ddd{
\lim_{x\to 0} \int_{k\in \Gamma} e^{ikx}  \, \frac{\cos(\omega t)}{k}  \,   dk
&=
\lim_{x\to 0} \int_{k\in \mathcal C} e^{ikx}  \, \frac{\cos(\omega t)}{k}  \,   dk
=
\int_{k\in \mathcal C}  \frac{\cos(\omega t)}{k}  \,   dk
=
i\pi \left( \cos t - 1\right).
}
Furthermore, noting that $C_\varepsilon$ is oriented clockwise, we compute
\eee{
\lim_{x\to 0} \int_{k\in C_\varepsilon} e^{ikx}  \, \frac{\cos(\omega t)}{k}  \,   dk
=
\int_{k\in C_\varepsilon}  \frac{\cos(\omega t)}{k}  \,   dk
=
-i\pi \text{Res}\left[\frac{\cos(\omega t)}{k}, k=0\right]
=
-i\pi.
}
Thus, 
\eee{
\lim_{x\to 0} \dashint_{k\in \mathbb R} e^{ikx}  \, \frac{\cos(\omega t)}{k}  \,   dk
=
i\pi  \cos t
}
and, in  turn, \eqref{r-eq-sol-ver3} becomes
\eee{\label{r-eq-sol-ver4}
q(0, t)
=
\frac{1}{2i\pi}  \dashint_{k\in \mathbb R}   \frac{\cos(\omega t)}{k}  \,   \what{q_0'}(k)  dk
+
\frac{1}{2i\pi}   \int_{k\in \mathcal C}  \frac{\cos(\omega t)}{k}  \, \what{q_0'}(-k)   dk
-\frac 12 g_0(0) + g_0(t),
}
where we have passed the limit $x\to 0$ inside the remaining principal value integral since this integral converges absolutely (i.e. without the help from $e^{ikx}$) thanks to the fact that, as $|k|\to \infty$, $\omega \simeq 1$ (and hence $\cos(\omega t)$ is bounded) and $\what{q_0'}(k) =  \frac{1}{ik}q_0'(0) + \frac{1}{ik} \what{q_0''}(k)$ (recall that $q_0$ belongs to the Schwartz class).

Finally, invoking Cauchy's theorem and Jordan's lemma once again,\footnote{The exponential decay required for Jordan's lemma is provided by the half-line Fourier transform $\what{r_0'}(-k) := \int_{x=0}^\infty e^{ikx} r_0'(x) dx$.} we have
\eee{
\int_{k\in \mathcal C}  \frac{\cos(\omega t)}{k}  \, \what{q_0'}(-k)   dk
=
\dashint_{k\in \mathbb R}  \frac{\cos(\omega t)}{k}  \, \what{q_0'}(-k)   dk + \int_{k\in C_\varepsilon}  \frac{\cos(\omega t)}{k}  \, \what{q_0'}(-k)   dk.
}
Therefore, using the change of variable $k\mapsto -k$ as appropriate, we find
\ddd{\label{r-eq-sol-ver5}
q(0, t)
&=
\frac{1}{2i\pi}  \dashint_{k\in \mathbb R}   \frac{\cos(\omega t)}{k}  \,   \what{q_0'}(k)  dk
+
\frac{1}{2i\pi}  \dashint_{k\in \mathbb R}  \frac{\cos(\omega t)}{k}  \, \what{q_0'}(-k)   dk
\nn\\
&\quad
 + \frac{1}{2i\pi}  \int_{k\in C_\varepsilon}  \frac{\cos(\omega t)}{k}  \, \what{q_0'}(-k)   dk
-\frac 12 g_0(0) + g_0(t)
\nn\\
&=
\frac{1}{2i\pi}  \dashint_{k\in \mathbb R}   \frac{\cos(\omega t)}{k}  \,   \what{q_0'}(k)  dk
+
\frac{1}{2i\pi}  \dashint_{k\in \mathbb R}  \frac{\cos(\omega t)}{-k}  \, \what{q_0'}(k)   dk
\nn\\
&\quad
 + \frac{1}{2i\pi}  \int_{k\in C_\varepsilon}  \frac{\cos(\omega t)}{k}  \, \what{q_0'}(-k)   dk
-\frac 12 g_0(0) + g_0(t)
\nn\\
&=
\frac{1}{2i\pi}  \int_{k\in C_\varepsilon}  \frac{\cos(\omega t)}{k}  \, \what{q_0'}(-k)   dk
-\frac 12 g_0(0) + g_0(t).
}
Finally, by Cauchy's residue theorem we compute
\ddd{
\frac{1}{2i\pi}  \int_{k\in C_\varepsilon}  \frac{\cos(\omega t)}{k}  \, \what{q_0'}(-k)   dk
&=
-\frac 12 \text{Res}\left[\frac{\cos(\omega t)}{k}  \, \what{q_0'}(-k), k=0\right]
\nn\\
&=
- \frac 12 \what{q_0'}(0)
:=
-\frac 12 \int_{x=0}^\infty q_0'(x) dx
= \frac 12 q_0(0)
= \frac 12 g_0(0)
}
with the last equality due to the compatibility condition $q_0(0) = q(0, 0) = g_0(0)$.
Hence, we overall conclude that if $q(x, t)$ is defined by the UTM formula \eqref{lcbous-sol-fin-q} then $q(0, t) = g_0(t)$, i.e. the boundary condition \eqref{lbous-bc} is satisfied.

\begin{remark}[Uniform convergence at the boundary]
As a final remark, we emphasize the importance of the complex contour $\mathcal C$ in the UTM solution formulae \eqref{lcbous-sol-fin}. Indeed, as illustrated by the above computations, had we deformed from $\mathcal C$ back to $\mathbb R$ it would not have been possible to explicitly verify that our formulae satisfy IBVP \eqref{lbous-ibvp} (and, in particular, the boundary condition) due to the loss of uniform convergence at the boundary induced by the deformation to the real axis.
\end{remark}

\section{Conclusion}
\label{conc-s}

A novel solution formula for the linearized classical Boussinesq system on the half-line with Dirichlet boundary data was derived by employing the unified transform method of Fokas. 
More precisely, the analysis utilized the recently formulated extension of the linear component of Fokas's method from single equations to systems \cite{dgsv2018}, as well as fundamental ideas of the method such as escaping to the complex spectral plane and exploiting the symmetries of a central identity known as the global relation. 

\begin{figure}[ht!]
\begin{center}
\includegraphics[height=30mm,width=49mm]{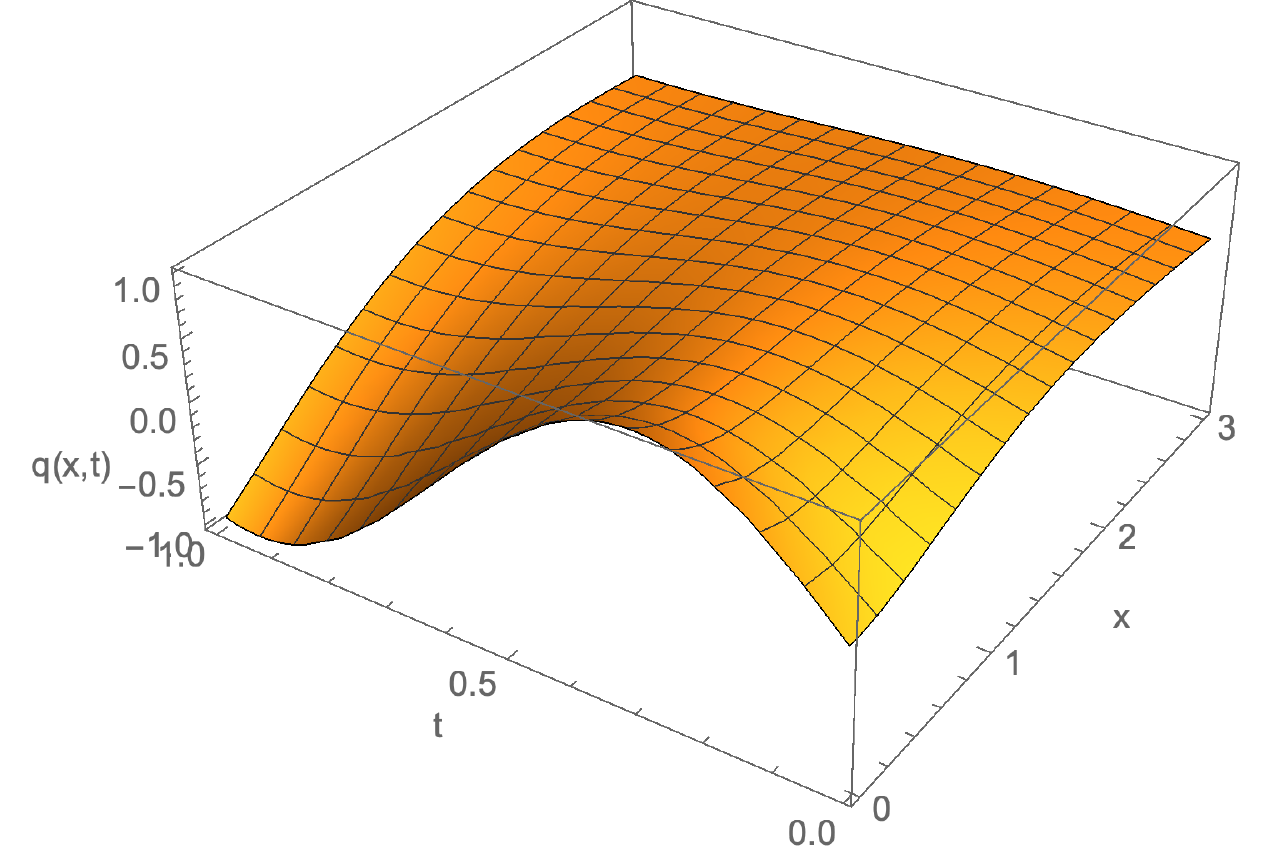}
\hspace*{3mm}
\includegraphics[height=29mm,width=47mm]{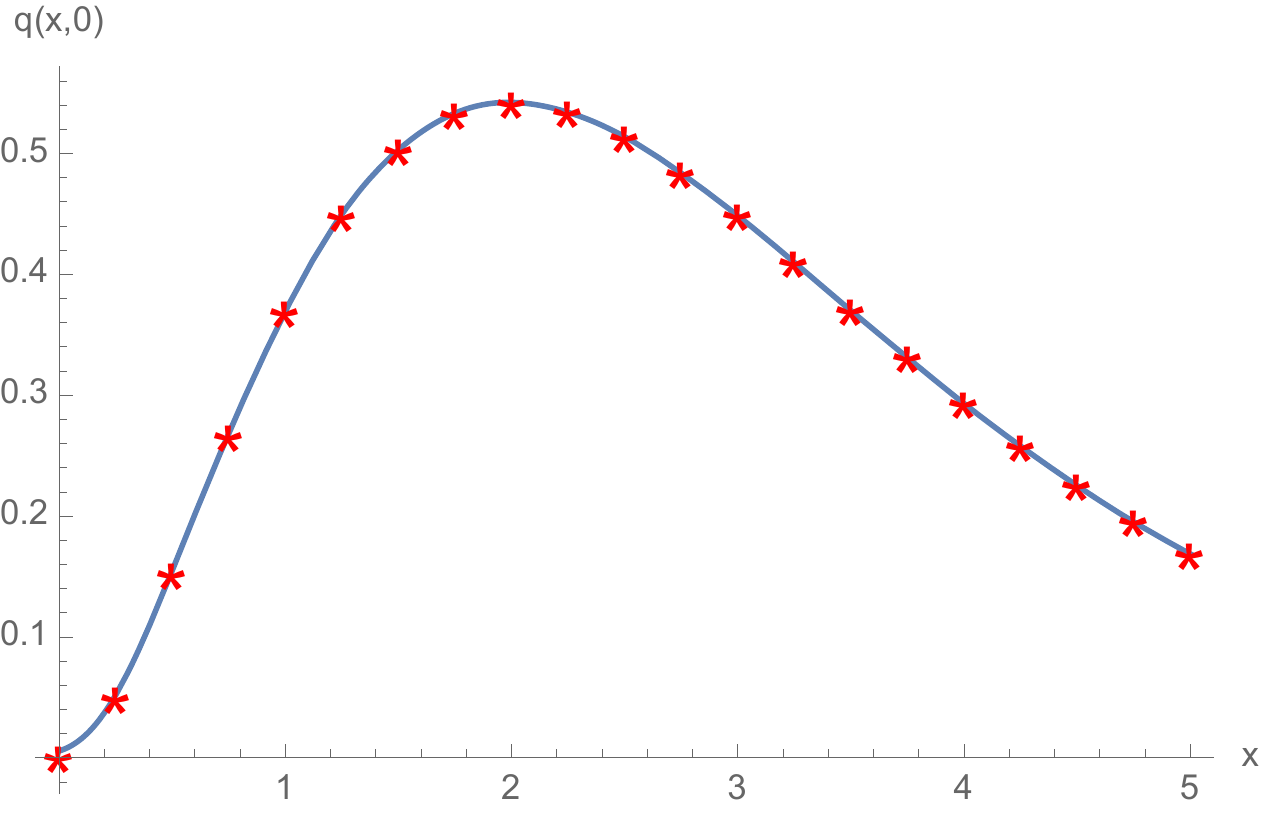}
\hspace*{3mm}
\includegraphics[height=29mm,width=47mm]{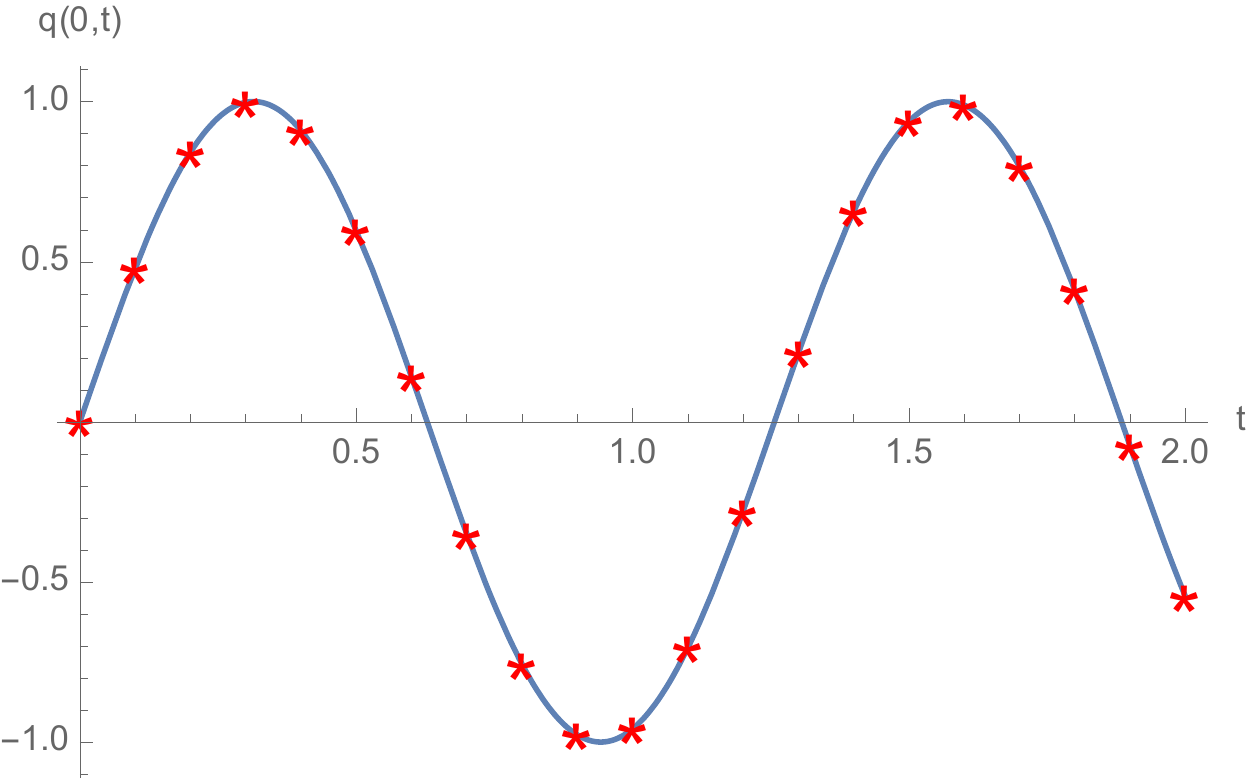}
\end{center}
\caption{Evaluation of the novel solution formula \eqref{lcbous-sol-fin-q} for the $q$-component of IBVP \eqref{lbous-ibvp} with initial data $q_0(x) = x^2 e^{-x}$ and boundary data $g_0(t) = \sin(5t)$. 
\textit{Left~panel:} The solution $q$ as a function of $(x, t)\in (0, 3) \times (0, 1)$.
\textit{Center panel:} The solution $q$ at $t=0$ as a function of $x \in (0, 5)$.
\textit{Right panel:} The solution $q$ at $x=0$ as a function of $t \in (0, 2)$.
The red stars in the center and right panels represent the \textit{precise} values of the data $q_0(x)$ and $g_0(t)$, respectively, showing perfect agreement of formula \eqref{lcbous-sol-fin-q} with the prescribed data. Of particular importance for the third plot is the uniform convergence of formula \eqref{lcbous-sol-fin-q} at the boundary $x=0$.}
\label{plots-f}
\end{figure}

The resulting solution formula \eqref{lcbous-sol-fin} has the important advantage of uniform convergence at the boundary $x=0$, thereby allowing for its explicit verification against the linearized classical Boussinesq IBVP \eqref{lbous-ibvp} through a direct calculation (see Section \ref{ver-s}).  
Beyond uniform convergence, the novel formula enjoys exponentially decaying integrands and hence, as usual with formulae derived via the unified transform, it is expected to be particularly effective regarding numerical considerations. Indeed, only a few lines of code in Mathematica result in the three plots of Figure~\ref{plots-f}.

Furthermore, formula \eqref{lcbous-sol-fin-q} was rederived in Section \ref{one-eq-s} from a different starting point, namely by reducing system \eqref{lbous-sys} to a single equation. This reduction is not possible for general dispersive systems, which is the main reason why the systems approach of Section \ref{gr-s} is preferable; however, the second approach demonstrates the versatility of the unified transform method and offers a different perspective concerning the types of admissible boundary data. Moreover, to the best of our knowledge, the analysis of Section \ref{one-eq-s} signifies the first time that the unified transform method has been applied to an equation whose dispersion relation is a quotient involving a complex square root (and hence branching), due to the presence of the mixed derivative  $\left(1-\p_x^2\right)\p_t^2$.

We emphasize that the present article is not the first one devoted to the linearized classical Boussinesq IBVP \eqref{lbous-ibvp} via the unified transform method. Indeed, Fokas and Pelloni  had previously considered the same problem in \cite{fp2005}. Importantly, however, their approach was entirely different, as they utilized the \textit{nonlinear} component of the unified transform method, which relies on formulating IBVP \eqref{lbous-ibvp} as a Lax pair and then integrating it by using ideas inspired from the inverse scattering transform, namely by associating it with a Riemann-Hilbert problem which is then solved via Plemelj's formulae. The complexity of those techniques naturally limits the accessibility of the derivation of \cite{fp2005} to a very specialized audience. In contrast, the present work employs the \textit{linear} component of the unified transform method, which only requires knowledge of Fourier transform and Cauchy's theorem from complex analysis and hence is accessible to a broad audience within the applied sciences.

Another significant difference between the present work and \cite{fp2005} is the fact that, due to the nature of the Riemann-Hilbert problem, the resulting solution formula in \cite{fp2005} involves certain principal value integrals. The presence of these integrals seems unsuitable for the purpose of  effective numerical implementations. Furthermore, it will most likely pose an issue when attempting to use that formula for establishing well-posedness of the original, nonlinear classical Boussinesq system \eqref{cbous} via the contraction mapping approach (see relevant discussion in the introduction). On the other hand, the novel solution formula derived in this article does not involve any singularities and hence is expected to be effective for showing well-posedness of the classical Boussinesq system \eqref{cbous} on the half-line via the contraction mapping approach, along the lines of \cite{fhm2016,fhm2017,hm2015b,hm2020}.

\vskip 4mm

\noindent
\textbf{Acknowledgements.} The first  author would like to thank the Department of Mathematics of the University of Kansas for partially supporting their research  through an undergraduate research award. All three authors are grateful to Andre Kurait for inspiring discussions during the 2018-19 academic year that paved the way to the present work.  Finally, the authors are thankful to the  reviewers of the manuscript for useful remarks and suggestions.

\end{document}